\input amstex
\documentstyle{amsppt}
%
\catcode`@=11
\redefine\output@{%
  \def\break{\penalty-\@M}\let\par\endgraf
  \ifodd\pageno\global\hoffset=105pt\else\global\hoffset=8pt\fi  
  \shipout\vbox{%
    \ifplain@
      \let\makeheadline\relax \let\makefootline\relax
    \else
      \iffirstpage@ \global\firstpage@false
        \let\rightheadline\frheadline
        \let\leftheadline\flheadline
      \else
        \ifrunheads@ 
        \else \let\makeheadline\relax
        \fi
      \fi
    \fi
    \makeheadline \pagebody \makefootline}%
  \advancepageno \ifnum\outputpenalty>-\@MM\else\dosupereject\fi
}
\catcode`\@=\active
\nopagenumbers
\def\negskp{\hskip -2pt}
\def\Alpha{\operatorname{A}}
 \accentedsymbol\bundleboldX{\raise 8pt \vbox{\hsize=7pt\noindent
  \vrule height 0.4pt depth 3pt width 1pt\vbox{\hrule width 5pt}}
  \kern -5pt\bold X}
\accentedsymbol\bd{\kern 2pt\bar{\kern -2pt d}}
\accentedsymbol\bbd{\kern 2pt\bar{\kern -2pt\bold d}}

\def\blue#1{#1}
\def\red#1{#1}
\catcode`#=11\def\diez{#}\catcode`#=6
\catcode`_=11\def\podcherkivanie{_}\catcode`_=8
\def\mycite#1{\cite{\blue{#1}}\immediate\special{ps:
     ShrHPSdict begin /ShrBORDERthickness 0 def}}

\def\mytag#1{%
    \tag#1}
\def\mythetag#1{\thetag{\blue{#1}}\immediate\special{ps:
     ShrHPSdict begin /ShrBORDERthickness 0 def}}
\def\myrefno#1{\no#1}
\def\myhref#1#2{\blue{#2}\immediate\special{ps:
     ShrHPSdict begin /ShrBORDERthickness 0 def}}
\def\myEarXivlink{\myhref{http://arXiv.org}{http:/\negskp/arXiv.org}}
\def\myGeoCities{\myhref{http://www.geocities.com}{GeoCities}}
\def\mytheorem#1{\csname proclaim\endcsname{Theorem #1}}

\def\mylemma#1{\csname proclaim\endcsname{Lemma #1}}

\def\mycorollary#1{\csname proclaim\endcsname{Corollary #1}}

\def\mydefinition#1{\definition{Definition #1}}
\def\mythedefinition#1{\blue{#1}\immediate\special{ps:
     ShrHPSdict begin /ShrBORDERthickness 0 def}}

\pagewidth{360pt}
\pageheight{606pt}
\topmatter
\title
A cubic identity for the Infeld-van der Waerden field
and its application.
\endtitle
\author
R.~A.~Sharipov
\endauthor
\address 5 Rabochaya street, 450003 Ufa, Russia\newline
\vphantom{a}\kern 12pt Cell Phone: +7(917)476 93 48
\endaddress
\email \vtop to 30pt{\hsize=280pt\noindent
\myhref{mailto:r-sharipov\@mail.ru}
{r-sharipov\@mail.ru}\newline
\myhref{mailto:R\podcherkivanie Sharipov\@ic.bashedu.ru}
{R\_\hskip 1pt Sharipov\@ic.bashedu.ru}\vss}
\endemail
\urladdr
\vtop to 20pt{\hsize=280pt\noindent
\myhref{http://www.geocities.com/r-sharipov}
{http:/\negskp/www.geocities.com/r-sharipov}\newline
\myhref{http://www.freetextbooks.boom.ru/index.html}
{http:/\negskp/www.freetextbooks.boom.ru/index.html}\vss}
\endurladdr
\abstract
    A cubic identity for the Infeld-van der Waerden field is found
and its application to verifying an explicit formula for the spinor
components of the metric connection is demonstrated.
\endabstract
\endtopmatter
\loadbold
\TagsOnRight
\document
\accentedsymbol\hatboldsymbolgamma{\kern 1.3pt\hat{\kern -1.3pt
   \boldsymbol\gamma}}
\accentedsymbol\hatgamma{\kern 1.3pt\hat{\kern -1.3pt\gamma}}

\rightheadtext{A cubic identity for the Infeld-van der Waerden field \dots}
\head
1. Introduction. 
\endhead
    The Infeld-van der Waerden field $\bold G$ is a special spin-tensorial 
field associated with the bundle of Weyl spinors $SM$ over the four-dimensional 
space-time manifold $M$ in general relativity. Along with the spinor metric 
$\bold d$, the Infeld-van der Waerden field $\bold G$ forms the basic equipment
of the spinor bundle $SM$. In this role $\bold d$ and $\bold G$ are similar to
the metric tensor $\bold g$ of $M$.
    The fields $\bold d$, $\bold G$, and $\bold g$ are related to each other
through a series of identities, e\.\,g. we have
$$
\xalignat 2
&\hskip -2em
\sum^2_{r=1}\sum^2_{\bar r=1}G^{r\kern 0.5pt\bar r}_p\ 
G^{\,q}_{r\kern 0.5pt\bar r}=2\,\delta^q_p,
&&\sum^3_{q=0}G^{r\kern 0.5pt\bar r}_q\ 
G^{\,q}_{s\kern 0.5pt\bar s}=2\,\delta^r_s
\,\delta^{\bar r}_{\bar s}.\qquad
\mytag{1.1}
\endxalignat
$$
The identities \mythetag{1.1} were considered in \mycite{1}. They 
are quadratic with respect to the components of the Infeld-van der Waerden 
field $\bold G$. Here we consider the identity
$$
\hskip -2em
\gathered
\sum^2_{s=1}\sum^2_{\bar s=1}G^{r\kern 0.5pt\bar s}_p\,
G^{\,m}_{s\kern 0.5pt\bar s}\,G^{s\kern 0.5pt\bar r}_q
=G^{r\kern 0.5pt\bar r}_p\,\delta^m_q+G^{r\kern 0.5pt\bar r}_q
\,\delta^m_p-\sum^3_{n=0}G^{r\kern 0.5pt\bar r}_n\,g^{mn}
\,g_{p\kern 1pt q}\,+\\
+\sum^3_{a=0}\sum^3_{b=0}\sum^3_{n=0}i\,g_{pa}\,g_{qb}
\,\omega^{ambn}\,G^{r\kern 0.5pt\bar r}_n,
\endgathered
\mytag{1.2}
$$
which is cubic with respect to the components of $\bold G$.\par
     The metric tensor $\bold g$ is canonically associated with its metric
connection $\Gamma$. The metric connection $\Gamma$, in its turn, has an 
extension $(\Gamma,\Alpha,\bar{\Alpha})$ to the spinor bundle $SM$. The 
spinor components $\Alpha^i_{\kern 0.5pt r j}$ of this extension are given by 
the formula
$$
\pagebreak
\hskip -2em
\gathered
\Alpha^i_{\kern 0.5pt r j}
=\sum^2_{\bar s=1}\sum^3_{p=0}\sum^3_{q=0}
\frac{G^{i\bar s}_p\,\Gamma^p_{rq}
\,G^{\,q}_{j\kern 0.5pt\bar s}}{4}\,-\\
-\sum^2_{\bar s=1}\sum^3_{q=0}
\frac{L_{\boldsymbol\Upsilon_{\!r}}\!(G^{i\bar s}_q)
\,G^{\,q}_{j\kern 0.5pt\bar s}}{4}
-\sum^2_{\bar i=1}\sum^2_{\bar j=1}
\frac{L_{\boldsymbol\Upsilon_{\!r}}\!(\bd_{\bar j\kern 0.5pt\bar i})
\,\bd^{\kern 0.5pt\bar i\bar j}\,\delta^{\,i}_j}{4}.
\endgathered
\mytag{1.3}
$$
The formula \mythetag{1.3} was derived in \mycite{1} by means of direct,
but rather huge calculations. In this paper we verify this formula with the
use of the identity \mythetag{1.2}.
\head
2. Coordinates and frames.
\endhead
     The components of a vector and the components of a tensor are always 
relative to some basis. In the case of vectorial and tensorial fields we
should have bases at each point of the space-time manifold $M$. Thus we
come to the concept of a frame.
\mydefinition{1.1} A frame of the tangent bundle $TM$ is a quadruple of
vector fields $\boldsymbol\Upsilon_0,\,\boldsymbol\Upsilon_1,\,\boldsymbol
\Upsilon_2,\,\boldsymbol\Upsilon_3$ defined in some open domain of $M$ and
linearly independent at each point of this domain. 
\enddefinition
     Let $\boldsymbol\Upsilon_0,\,\boldsymbol\Upsilon_1,\,\boldsymbol
\Upsilon_2,\,\boldsymbol\Upsilon_3$ be a frame of the tangent bundle $TM$
Then the commutator $[\boldsymbol\Upsilon_{\!i},\boldsymbol\Upsilon_{\!j}]$
is a vector field that can be expressed through the frame fields:
$$
\hskip -2em
[\boldsymbol\Upsilon_{\!i},\boldsymbol\Upsilon_{\!j}]=\sum^3_{k=0}
c^{\,k}_{ij}\,\boldsymbol\Upsilon_{\!k}.
\mytag{2.1}
$$
The coefficients $c^{\,k}_{ij}$ in \mythetag{2.1} are called the 
{\it commutation coefficients\/} of the frame $\boldsymbol\Upsilon_0,
\,\boldsymbol\Upsilon_1$, $\boldsymbol\Upsilon_2,\,\boldsymbol\Upsilon_3$. 
They are uniquely determined for any frame.
\mydefinition{1.2} A frame $\boldsymbol\Upsilon_0,\,\boldsymbol
\Upsilon_1,\,\boldsymbol\Upsilon_2,\,\boldsymbol\Upsilon_3$ of the 
tangent bundle $TM$ is a called a {\it holonomic frame\/} if all of
its commutation coefficients $c^{\,k}_{ij}$ are identically zero. 
\enddefinition
      Once we take some local coordinates $x^0,\,x^1,\,x^2,\,x^3$ in
the space-time manifold $M$, we get the holonomic frame of the coordinate 
vector fields 
$$
\xalignat 4
&\hskip -2em
\boldsymbol\Upsilon_0=\frac{\partial}{\partial x^0},
&&\boldsymbol\Upsilon_1=\frac{\partial}{\partial x^1},
&&\boldsymbol\Upsilon_2=\frac{\partial}{\partial x^2},
&&\boldsymbol\Upsilon_3=\frac{\partial}{\partial x^3}
\qquad
\mytag{2.2}
\endxalignat
$$
defined within the domain of the coordinates $x^0,\,x^1,\,x^2,\,x^3$ 
and naturally associated to them. Local coordinates and their coordinate
frames \mythetag{2.2} are sufficient for most calculations in general
relativity where no spinors are involved. However, when dealing with
spinors we need to use non-holonomic frames.
\mydefinition{1.3} A frame of the spinor bundle $SM$ is a pair of
smooth sections $\boldsymbol\Psi_1,\,\boldsymbol\Psi_2$ of $SM$ defined 
in some open domain of $M$ and linearly independent at each point of 
this domain. 
\enddefinition
     Frames of $TM$ are often called {\it spacial frames\/}, while 
frames of $SM$ are called {\it spinor frames}. Having some spacial frame 
$\boldsymbol\Upsilon_0,\,\boldsymbol\Upsilon_1,\,\boldsymbol\Upsilon_2,
\,\boldsymbol\Upsilon_3$ and some spinor frame $\boldsymbol\Psi_1,\,
\boldsymbol\Psi_2$ with common domain, we can express tensorial and spin
tensorial fields through their components forming multi-indexed arrays.
The indices of such arrays are subdivided into three groups: 
{\it spacial indices\/} ranging from $0$ to $3$, {\it spinor indices},
and {\it conjugate spinor indices\/} both ranging from $1$ to $2$. 
Indices within each group can be either upper or lower. The quantities 
$g_{p\kern 1pt q}$ with two lower spacial indices in \mythetag{1.2} are
the components of the metric tensor $\bold g$, while the quantities 
$g^{mn}$ are the components of the dual metric tensor. The components
of these two tensors form two $4\times 4$ symmetric matrices inverse 
to each other:
$$
\hskip -2em
\sum^3_{r=0}g_{p\kern 1pt r}\,g^{r\kern 0.5pt q}=\delta^{\kern 0.5ptq}_p.
\mytag{2.3}
$$
Similarly, $d_{ij}$ and $d^{\kern 0.5pt ij}$ are the components of the 
spinor metric and dual spinor metric respectively. They form $2\times 2$ 
skew-symmetric matrices inverse to each other:
$$
\hskip -2em
\sum^2_{r=1}d_{jr}\,d^{\kern 0.5pt ri}=\delta^i_j.
\mytag{2.4}
$$\par
     The Infeld-van der Waerden field $\bold G$ in \mythetag{1.1} is
presented by the quantities $G^{r\kern 0.5pt\bar r}_p$. It has one upper
spinor index $r$, one upper conjugate spinor index $\bar r$, and one
spacial index $p$. The quantities $G^{\,q}_{s\kern 0.5pt\bar s}$ are the
components of the inverse Infeld-van der Waerden field. They are produced 
from $G^{r\kern 0.5pt\bar r}_p$ by raising the spacial index $p$ and
lowering the spinor and conjugate spinor indices $r$ and $\bar r$:
$$
\hskip -2em
G^{\,q}_{s\kern 0.5pt\bar s}=\sum^2_{r=1}\sum^2_{\bar r=1}\sum^3_{p=0}
G^{r\kern 0.5pt\bar r}_p\,g^{p\kern 1 pt q}\,d_{rs}\,\bd_{\bar r\bar s}.
\mytag{2.5}
$$
The quantities $\bd_{\bar r\bar s}$ in the formula \mythetag{2.5} are the 
components of the conjugate spinor metric $\bbd=\tau(\bold d)$. They are 
related to $d_{ij}$ by means of the complex conjugation:
$$
\hskip -2em
\bd_{\,\bar i\kern 0.5pt\bar j}=\overline{d_{\,\bar i\kern 0.5pt\bar j}}.
\mytag{2.6}
$$
Similarly, the quantities $\bd^{\kern 0.5pt\bar i\bar j}$ in the formula
\mythetag{1.3} are defined as follows:
$$
\hskip -2em
\bd^{\kern 0.5pt\bar i\bar j}=\overline{d^{\raise 1.2pt \hbox{$\ssize\kern 
0.5pt\bar i\bar j$}}}.
\mytag{2.7}
$$
The quantities \mythetag{2.6} and \mythetag{2.7} form two mutually inverse
matrices:
$$
\hskip -2em
\sum^2_{r=1}\bd_{\bar j\bar r}\,\bd^{\kern 0.5pt\bar r\bar i}=\delta^i_j.
\mytag{2.8}
$$
The equality \mythetag{2.8} is similar to the above equalities 
\mythetag{2.3} and \mythetag{2.4}.\par
     The quantities $G^{r\kern 0.5pt\bar r}_p$ and $G^{\,q}_{s\kern 0.5pt
\bar s}$ satisfy the following identities:
$$
\xalignat 2
\hskip -2em
&G^{r\kern 0.5pt\bar r}_p=\overline{G^{\raise 1.2pt \hbox{$\ssize\bar r
\kern 0.5pt r$}}_p},
&&G^{\,q}_{s\kern 0.5pt\bar s}=\overline{G^{\,q}_{\bar s\kern 0.5pt s}}.
\mytag{2.9}
\endxalignat
$$
The identities \mythetag{2.9} mean that for each fixed $p$ and for each fixed
$q$ the quantities $G^{r\kern 0.5pt\bar r}_p$ and $G^{\,q}_{s\kern 0.5pt\bar s}$
form two Hermitian matrices. In the coordinate free form the identities 
\mythetag{2.9} are written as $\tau(\bold G)=\bold G$.\par
     Now let's proceed to the quantities $\omega^{ambn}$ in the formula 
\mythetag{1.2}. They are the components of the {\it dual volume tensor}.
These quantities are produced from the components of the {\it volume tensor}
$\boldsymbol\omega$ by means of the index raising procedure:
$$
\hskip -2em
\omega^{ambn}=\sum^3_{p=0}\sum^3_{p=0}\sum^3_{r=0}\sum^3_{s=0}
\omega_{prqs}\,g^{pa}\,g^{rm}\,g^{qb}\,g^{sn}.
\mytag{2.10}
$$
Both quantities $\omega^{ambn}$ and $\omega_{prqs}$ in \mythetag{2.10} can
be produced from the components of the metric tensor and the components of 
the dual metric tensor with the use of the Levi-Civita symbol. Indeed, we
have the following formula for $\omega_{prqs}$:
$$
\hskip -2em
\omega_{prqs}=\sqrt{-\det(g_{ij})}\ \varepsilon_{prqs}.
\mytag{2.11}
$$
Similarly, for $\omega^{ambn}$ we have the formula
$$
\hskip -2em
\omega^{\kern 0.5pt ambn}=-\sqrt{-\det(g^{ij})}\ \varepsilon^{ambn}.
\mytag{2.12}
$$
Regardless of the upper and lower positions of the indices, the components
of the Levi-Civita symbol are given by the formula
$$
\varepsilon^{ambn}=\varepsilon_{ambn}=\cases\ \ 0 &\vtop{\offinterlineskip
\hsize =140pt\noindent if at least two of the four indices $a\,b\,m\,n$ 
do coincide;\strut}\\
\ \ 1 & \vtop{\offinterlineskip\hsize =155pt\noindent if the indices $a\,m\,
b\,n$ form an even transposition of the numbers $0\,1\,2\,3$;\strut}\\
-1 & \vtop{\offinterlineskip\hsize =155pt\noindent if the indices $a\,m\,
b\,n$ form an odd transposition of the numbers $0\,1\,2\,3$.\strut}\\
\endcases
\qquad
\mytag{2.13}
$$
Note that the space-time manifold $M$ is assumed to be equipped with 
three geometric structures: the {\bf metric}, the {\bf orientation}, 
and the {\bf polarization} (see \mycite{2}). The orientation in $M$
means that we can distinguish right and left frames. The formulas 
\mythetag{2.11} and \mythetag{2.12} are written for right frames. 
Passing to left frames, we should change the sign in the right hand 
sides of them.\par
     Let's consider the Lie derivatives $L_{\boldsymbol\Upsilon_{\!r}}$ in
the formula \mythetag{1.3}. In the case of a holonomic frame \mythetag{2.2} 
the operators $L_{\boldsymbol\Upsilon_{\!r}}$ coincide with the corresponding 
partial derivatives, i\.\,e\. $L_{\boldsymbol\Upsilon_{\!r}}=\partial/\partial 
x^r$. In the case of a non-holonomic frame we should expand $\boldsymbol
\Upsilon_0,\,\boldsymbol\Upsilon_1,\,\boldsymbol\Upsilon_2,\,\boldsymbol
\Upsilon_3$ in some auxiliary holonomic frame:
$$
\hskip -2em
\boldsymbol\Upsilon_{\!r}=\sum^3_{i=0}\Upsilon^i_{\!r}\,\frac{\partial}
{\partial x^i}.
\mytag{2.14}
$$
Then the operators $L_{\boldsymbol\Upsilon_{\!r}}$ in \mythetag{1.3} act
as the differential operators in the right hand sides of the formula 
\mythetag{2.14}. In other words, for $L_{\boldsymbol\Upsilon_{\!r}}$ we 
have the expression 
$$
\hskip -2em
L_{\boldsymbol\Upsilon_{\!r}}=\sum^3_{i=0}\Upsilon^i_{\!r}\,\frac{\partial}
{\partial x^i}
\mytag{2.15}
$$
formally coinciding with \mythetag{2.14}. The Lie derivatives \mythetag{2.15}
satisfy the relationships
$$
\hskip -2em
[L_{\boldsymbol\Upsilon_{\!i}},L_{\boldsymbol\Upsilon_{\!j}}]
=\sum^3_{k=0}c^{\,k}_{ij}\,L_{\boldsymbol\Upsilon_{\!k}}.
\mytag{2.16}
$$
The commutation coefficients $c^{\,k}_{ij}$ in \mythetag{2.16} are the same 
as in the formula \mythetag{2.1}.
\head
3. Proof of the identity \mythetag{1.2}.
\endhead
The proof of the identity \mythetag{1.2} is similar to that of the
identities \mythetag{1.1} in \mycite{1}. It is based on direct calculations.
The matter is that the spinor bundle $SM$ is related to the tangent bundle
$TM$ through canonically associated frame pairs, while the fields $\bold g$,
$\bold d$, $\bold G$, and $\boldsymbol\omega$ are given by explicit formulas 
in such frame pairs.
\mydefinition{3.1} A frame $\boldsymbol\Upsilon_0,\,\boldsymbol
\Upsilon_1,\,\boldsymbol\Upsilon_2,\,\boldsymbol\Upsilon_3$ of the 
tangent bundle $TM$ is a called an {\it orthonormal frame\/} if the 
metric tensor $\bold g$ and its dual metric tensor \pagebreak
are given by the standard Minkowski matrix in this frame:
$$
\hskip -2em
g_{ij}=g^{ij}=\Vmatrix 1 & 0 & 0 & 0\\0 & -1 & 0 & 0\\
0 & 0 & -1 & 0\\0 & 0 & 0 & -1\endVmatrix.
\mytag{3.1}
$$
\enddefinition
Let's recall again that $M$ is equipped with the {\bf metric}, the 
{\bf orientation}, and the {\bf polarization}. The polarization is 
a geometric structure that distinguishes the {\bf future half light 
cone} from the {\bf past half light cone}.
\mydefinition{3.2} A frame $\boldsymbol\Upsilon_0,\,\boldsymbol
\Upsilon_1,\,\boldsymbol\Upsilon_2,\,\boldsymbol\Upsilon_3$ of the 
tangent bundle $TM$ is a called {\it positively polarized\/} if 
its first vector field $\boldsymbol\Upsilon_0$ belongs to the 
interior of the future light half cone at each point of its domain. 
\enddefinition
\mydefinition{3.3} A frame $\boldsymbol\Psi_1,\,\boldsymbol\Psi_2$ 
of the spinor bundle $SM$ is a called {\it orthonormal} if the spinor 
metric $\bold d$ is given by the following matrix in this frame: 
$$
\hskip -2em
d_{ij}=\Vmatrix 0 & 1\\ 
\vspace{1ex} -1 & 0\endVmatrix.
\mytag{3.2}
$$
\enddefinition
\noindent
Due to the definition~\mythedefinition{3.3} the dual spinor metric 
is given by the matrix 
$$
d^{\kern 0.5pt ij}=\Vmatrix 0 & -1\\ 
\vspace{1ex} 1 & 0\endVmatrix
$$
in any orthonormal frame of the spinor bundle $SM$.\par
     According to the definition of the spinor bundle $SM$ (see \mycite{1}),
each positively polarized right orthonormal frame $\boldsymbol\Upsilon_0,
\,\boldsymbol\Upsilon_1,\,\boldsymbol\Upsilon_2,\,\boldsymbol\Upsilon_3$ 
of the tangent bundle $TM$ is associated with some orthonormal frame $\boldsymbol\Psi_1,\,\boldsymbol\Psi_2$ of the spinor bundle $SM$. Such 
frames form canonically associated frame pairs. In any canonically 
associated frame pair the Infeld-van der Waerden field $\bold G$ is 
given by the following Pauli matrices:
$$
\xalignat 2
&\hskip -2em
G^{i\kern 0.5pt\bar i}_0=\Vmatrix 1 & 0\\ 0 & 1\endVmatrix
=\sigma_0,
&&G^{i\kern 0.5pt\bar i}_2=\Vmatrix 0 & -i\\ i & 0\endVmatrix
=\sigma_2,\\
\vspace{-1.4ex}
&&&\mytag{3.3}\\
\vspace{-1.4ex}
&\hskip -2em
G^{i\kern 0.5pt\bar i}_1=\Vmatrix 0 & 1\\ 1 & 0\endVmatrix
=\sigma_1,
&&G^{i\kern 0.5pt\bar i}_3=\Vmatrix 1 & 0\\ 0 & -1\endVmatrix
=\sigma_3.
\endxalignat
$$
Substituting \mythetag{3.3}, \mythetag{3.2}, and \mythetag{3.1} into
the formula \mythetag{2.5}, we calculate the components of the inverse 
Infeld-van der Waerden field:
$$
\xalignat 2
&\hskip -2em
G^{\kern 0.5pt 0}_{i\kern 0.5pt\bar i}=\Vmatrix 1 & 0\\ 0 & 1\endVmatrix
=\sigma_0,
&&G^{\kern 0.5pt 2}_{i\kern 0.5pt\bar i}=\Vmatrix 0 & i\\ -i & 0\endVmatrix
=-\sigma_2,\\
\vspace{-1.4ex}
&&&\mytag{3.4}\\
\vspace{-1.4ex}
&\hskip -2em
G^{\kern 0.5pt 1}_{i\kern 0.5pt\bar i}=\Vmatrix 0 & 1\\ 1 & 0\endVmatrix
=\sigma_1,
&&G^{\kern 0.5pt 3}_{i\kern 0.5pt\bar i}=\Vmatrix 1 & 0\\ 0 & -1\endVmatrix
=\sigma_3.
\endxalignat
$$
And finally, substituting \mythetag{3.1} into \mythetag{2.12}, we calculate 
the components of the inverse volume tensor $\boldsymbol\omega$ in the right 
orthonormal frame $\boldsymbol\Upsilon_0,\,\boldsymbol\Upsilon_1,\,\boldsymbol
\Upsilon_2,\,\boldsymbol\Upsilon_3$:
$$
\hskip -2em
\omega^{\kern 0.5pt ambn}=-\varepsilon^{ambn}.
\mytag{3.5}
$$
The quantities $\varepsilon^{ambn}$ in \mythetag{3.5} are defined by the 
formula \mythetag{2.13}. Now in order to prove the identity \mythetag{1.2}
it is sufficient to substitute \mythetag{3.1}, \mythetag{3.3}, \mythetag{3.4}, 
and \mythetag{3.5} into \mythetag{1.2} and verify that this equality is 
valid for all particular values of the indices $p$, $q$, $m$, $r$, and $\bar r$. 
The following computer code does it for us:
\medskip
\leftline{\tt\red{\ Verification\podcherkivanie List:=[]:}}
\leftline{\tt\red{\ for p from 0 by 1 to 3 do}}
\leftline{\tt\red{\ \ for m from 0 by 1 to 3 do}}
\leftline{\tt\red{\ \ \ for q from 0 by 1 to 3 do}}
\leftline{\tt\red{\ \ \ \ for r from 1 by 1 to 2 do}}
\leftline{\tt\red{\ \ \ \ \ for br from 1 by 1 to 2 do}}
\leftline{\tt\red{\ \ \ \ \ Equ:=0:}}
\leftline{\tt\red{\ \ \ \ \ \ for s from 1 by 1 to 2 do}}
\leftline{\tt\red{\ \ \ \ \ \ \ for bs from 1 by 1 to 2 do}}
\leftline{\tt\red{\ \ \ \ \ \ \ \ Equ:=Equ+G\podcherkivanie ortho[p][r,bs]%
*Inv\podcherkivanie G\podcherkivanie ortho[m][s,bs]}}
\leftline{\tt\red{\ \ \ \ \ \ \ \ \ \ \ \ \ \ \ \ *G\podcherkivanie%
ortho[q][s,br]:}}
\leftline{\tt\red{\ \ \ \ \ \ od od:}}
\leftline{\tt\red{\ \ \ \ \ \ Equ:=Equ-G\podcherkivanie ortho[p][r,br]%
*delta[m,q]}}
\leftline{\tt\red{\ \ \ \ \ \ \ \ \ \ \ \ \ \ -G\podcherkivanie ortho[q][r,br]%
*delta[m,p]:}}
\leftline{\tt\red{\ \ \ \ \ \ for n from 0 by 1 to 3 do}}
\leftline{\tt\red{\ \ \ \ \ \ \ Equ:=Equ+G\podcherkivanie ortho[n][r,br]%
*Inv\podcherkivanie g\podcherkivanie ortho[m,n]%
*g\podcherkivanie ortho[p,q]:}}
\leftline{\tt\red{\ \ \ \ \ \ od:}}
\leftline{\tt\red{\ \ \ \ \ \ for a from 0 by 1 to 3 do}}
\leftline{\tt\red{\ \ \ \ \ \ \ for b from 0 by 1 to 3 do}}
\leftline{\tt\red{\ \ \ \ \ \ \ \ for n from 0 by 1 to 3 do}}
\leftline{\tt\red{\ \ \ \ \ \ \ \ \ Equ:=Equ-I*g\podcherkivanie ortho[p,a]%
*g\podcherkivanie ortho[q,b]*Inv\podcherkivanie omega[a,m,b,n]}}
\leftline{\tt\red{\ \ \ \ \ \ \ \ \ \ \ \ \ \ \ \ \ \ \ *G\podcherkivanie%
ortho[n][r,br]:}}
\leftline{\tt\red{\ \ \ \ \ \ od od od:}}
\leftline{\tt\red{\ \ \ \ \ \ Verification\podcherkivanie List:=%
[op(Verification\podcherkivanie List),evalb(Equ=0)]:}}
\leftline{\tt\red{\ od od od od od:}}
\medskip\noindent
Upon executing the above code it is sufficient to type
\medskip\noindent
\red{{\tt\vphantom{>} print(Verification\podcherkivanie List):}}
\medskip\noindent
and find that the identity \mythetag{1.2} is proved for any canonically 
associated pair of frames. Due to the tensorial nature of this identity,
being proved for some particular frame pair, it remains valid for 
arbitrary pairs of frames.\par
     Note that the above code is designed for the Maple\footnotemark\ 
package. However, it can be easily adapted for other symbolic computation
packages.\footnotetext{\ Maple is a trademark of Waterloo Maple Inc.}
\head
4. Other relationships derived from \mythetag{1.2}.
\endhead
     Let's multiply both sides of the identity \mythetag{1.2} by 
$G^{u\bar u}_m$ and sum it over the index $m$. As a result, applying 
the second identity \mythetag{1.1}, we get
$$
\hskip -2em
\gathered
2\,G^{r\kern 0.5pt\bar u}_p\,G^{u\kern 0.5pt\bar r}_q
=G^{r\kern 0.5pt\bar r}_p\,G^{u\bar u}_q+G^{r\kern 0.5pt\bar r}_q
\,G^{u\bar u}_p-2\,d^{ru}\,\bd^{\bar r\bar u}\,g_{p\kern 1pt q}\,+\\
+\sum^3_{a=0}\sum^3_{b=0}\sum^3_{m=0}\sum^3_{n=0}i\,g_{pa}\,g_{qb}
\,\omega^{ambn}\,G^{u\bar u}_m\,G^{r\kern 0.5pt\bar r}_n.
\endgathered
\mytag{4.1}
$$
Symmetrizing the equality \mythetag{4.1} with respect to $p$ and $q$,
we derive 
$$
\hskip -2em
G^{r\kern 0.5pt\bar u}_p\,G^{u\kern 0.5pt\bar r}_q
+G^{r\kern 0.5pt\bar u}_q\,G^{u\kern 0.5pt\bar r}_p
=G^{r\kern 0.5pt\bar r}_p\,G^{u\bar u}_q
+G^{r\kern 0.5pt\bar r}_q\,G^{u\bar u}_p
-2\,d^{\kern 0.5pt ru}\,\bd^{\kern 0.5pt \bar r\bar u}\,g_{p\kern 1pt q}.
\mytag{4.2}
$$
Alternating the equality \mythetag{4.1} with respect to the same pair of
indices, we obtain
$$
\hskip -2em
G^{r\kern 0.5pt\bar u}_p\,G^{u\kern 0.5pt\bar r}_q
-G^{r\kern 0.5pt\bar u}_q\,G^{u\kern 0.5pt\bar r}_p
=\sum^3_{a=0}\sum^3_{b=0}\sum^3_{m=0}\sum^3_{n=0}i\,g_{pa}\,g_{qb}
\,\omega^{ambn}\,G^{u\bar u}_m\,G^{r\kern 0.5pt\bar r}_n.
\mytag{4.3}
$$
The identities \mythetag{4.2} and \mythetag{4.3} taken together are
equivalent to \mythetag{4.1}.
\head
5. The metric connection and its spinor components.
\endhead
     The metric connection $\Gamma$ in $M$ associated with the metric tensor
$\bold g$ is a torsion-free connection satisfying the condition
$$
\hskip -2em
\nabla\bold g=0.
\mytag{5.1}
$$
The equality \mythetag{5.1} is known as the concordance condition for the 
metric and connection. In the coordinate form the torsion-free condition 
is written as
$$
\hskip -2em
\Gamma^k_{ij}-\Gamma^k_{j\kern 0.5pt i}=c^{\,k}_{ij}.
\mytag{5.2}
$$
As for the concordance condition \mythetag{5.1}, it expands to
$$
\hskip -2em
L_{\boldsymbol\Upsilon_{\!r}}(g_{ij})-\sum^3_{k=0}
\Gamma^k_{ri}\,g_{kj}-\sum^3_{k=0}\Gamma^k_{rj}\,g_{ik}=0.
\mytag{5.3}
$$
The equations \mythetag{5.2} and \mythetag{5.3} can be solved with respect
to $\Gamma^k_{ij}$. They yield
$$
\hskip -2em
\gathered
\Gamma^k_{ij}=\sum^3_{r=0}\frac{g^{\kern 0.5pt kr}}{2}
\left(L_{\boldsymbol\Upsilon_{\!i}}\!(g_{rj})
+L_{\boldsymbol\Upsilon_{\!j}}\!(g_{i\kern 0.5pt r})
-L_{\boldsymbol\Upsilon_{\!r}}\!(g_{ij})\right)+\\
+\,\frac{c^{\,k}_{ij}}{2}
-\sum^3_{r=0}\sum^3_{s=0}\frac{c^{\,s}_{i\kern 0.5pt r}}{2}\,g^{kr}
\,g_{sj}-\sum^3_{r=0}\sum^3_{s=0}\frac{c^{\,s}_{j\kern 0.5ptr}}{2}
\,g^{kr}\,g_{s\kern 0.5pt i}.
\endgathered
\mytag{5.4}
$$
The only term, which is skew-symmetric with respect to $i$ and $j$, is 
$c^{\,k}_{ij}/2$. Other terms are either symmetric or have their symmetric
counterparts. For this reason, substituting \mythetag{5.4} back into
\mythetag{5.2}, we find that this equality is identically fulfilled:
$$
\Gamma^k_{ij}-\Gamma^k_{j\kern 0.5pt i}=\frac{c^{\,k}_{ij}}{2}
-\frac{c^{\,k}_{j\kern 0.5pt i}}{2}=c^{\,k}_{ij}.
$$
Applying \mythetag{5.4} to the second term in \mythetag{5.3}, we derive
$$
\hskip -2em
\gathered
-\sum^3_{k=0}\Gamma^k_{ri}\,g_{kj}=-\frac{1}{2}
\left(L_{\boldsymbol\Upsilon_{\!r}}\!(g_{j\kern 0.5pt i})
+L_{\boldsymbol\Upsilon_{\!i}}\!(g_{r\kern -0.5pt j})
-L_{\boldsymbol\Upsilon_{\!j}}\!(g_{ri})\right)\,-\\
-\frac{1}{2}\sum^3_{k=0}c^{\,k}_{ri}\ g_{kj}
+\frac{1}{2}\sum^3_{s=0}c^{\,s}_{r\kern -0.5pt j}
\ g_{s\kern 0.5pt i}
+\frac{1}{2}\sum^3_{s=0}c^{\,s}_{ij}\ g_{s\kern 0.5pt r}.
\endgathered
\mytag{5.5}
$$
The third term in \mythetag{5.3} differs from the second term by exchanging 
$i$ and $j$:
$$
\hskip -2em
\gathered
-\sum^3_{k=0}\Gamma^k_{r\kern -0.5pt j}\,g_{ki}
=-\frac{1}{2}\left(L_{\boldsymbol\Upsilon_{\!r}}\!(g_{ij})
+L_{\boldsymbol\Upsilon_{\!j}}\!(g_{ri})
-L_{\boldsymbol\Upsilon_{\!i}}\!(g_{r\kern -0.5ptj})\right)\,-\\
-\frac{1}{2}\sum^3_{k=0}c^{\,k}_{r\kern -0.5ptj}\ g_{ki}
+\frac{1}{2}\sum^3_{s=0}c^{\,s}_{r\kern -0.5pt i}
\ g_{sj}+\frac{1}{2}\sum^3_{s=0}c^{\,s}_{j\kern 0.5pt i}
\ g_{s\kern 0.5pt r}.
\endgathered
\mytag{5.6}
$$
Adding the formulas \mythetag{5.5} and \mythetag{5.6}, we find that
$$
-\sum^3_{k=0}\Gamma^k_{ri}\,g_{kj}-\sum^3_{k=0}\Gamma^k_{r\kern -0.5pt j}
\,g_{ki}=-L_{\boldsymbol\Upsilon_{\!r}}\!(g_{ij}).
$$
This relationship is equivalent to \mythetag{5.3}. Thus, for the connection
with the components \mythetag{5.4} we have verified both conditions --- the 
torsion-free condition \mythetag{5.2} and the concordance condition 
\mythetag{5.3}.\par
     The metric connection $\Gamma$ with the components \mythetag{5.4} has
a spinor extension $(\Gamma,\Alpha,\bar{\Alpha})$ which is uniquely fixed
by the following two concordance conditions:
$$
\xalignat 2
&\hskip -2em
\nabla\bold d=0,
&&\nabla\bold G=0.
\mytag{5.7}
\endxalignat
$$
In the coordinate form the conditions \mythetag{5.7} are written as follows:
$$
\align
&\hskip -2em
L_{\boldsymbol\Upsilon_{\!r}}(d_{ij})-\sum^2_{k=1}
\Alpha^k_{ri}\,d_{kj}-\sum^2_{k=1}\Alpha^k_{rj}\,d_{i\kern 0.5pt k}=0,
\mytag{5.8}\\
&\hskip -2em
L_{\boldsymbol\Upsilon_{\!r}}(G^{i\kern 0.5pt\bar i}_p)
+\sum^2_{k=1}\Alpha^i_{rk}\,G^{k\kern 0.5pt \bar i}_p
+\sum^2_{\bar k=1}\overline{\Alpha^{\bar i}_{r\bar k}}
\,G^{i\kern 0.5pt \bar k}_p\,
-\sum^3_{k=0}\Gamma^k_{rp}\,G^{i\kern 0.5pt\bar i}_k=0.
\qquad
\mytag{5.9}
\endalign
$$
The $\Alpha$-components of the spinor extension of $\Gamma$ are given by the
formula \mythetag{1.3}. Our next goal is to verify the formula \mythetag{1.3}
by substituting it into \mythetag{5.8} and \mythetag{5.9}.\par
    Let's begin with the second term in the formula \mythetag{5.8}. Applying 
the formula \mythetag{1.3} to this term, we derive the following expression 
for it:
$$
\hskip -2em
\gathered
\sum^2_{k=1}\Alpha^k_{ri}\,d_{kj}=\frac{1}{4}\sum^2_{k=1}\sum^2_{\bar s=1}
\sum^3_{p=0}\sum^3_{q=0}G^{k\bar s}_p\,\Gamma^p_{rq}\,G^{\,q}_{i\kern 0.5pt
\bar s}\,d_{kj}\,-\\
-\,\frac{1}{4}\sum^2_{k=1}\sum^2_{\bar s=1}\sum^3_{q=0}L_{\boldsymbol
\Upsilon_{\!r}}\!(G^{k\bar s}_q)\,G^{\,q}_{i\kern 0.5pt\bar s}\,d_{kj}
-\frac{1}{4}\sum^2_{\bar i=1}\sum^2_{\bar j=1}L_{\boldsymbol\Upsilon_{\!r}}
\!(\bd_{\bar j\kern 0.5pt\bar i})\,\bd^{\kern 0.5pt\bar i\bar j}\,d_{ij}.
\endgathered
\mytag{5.10}
$$
For to transform the first term in the right hand side of \mythetag{5.10} 
we use the formulas
$$
\xalignat 2
&\hskip -2em
\sum^2_{k=1}G^{k\bar s}_p\,d_{kj}=\sum^3_{m=0}\sum^2_{\bar r=1}
G^{\,m}_{\!j\kern 0.5pt\bar r}\,\bd^{\kern 0.5pt\bar r\bar s}\,g_{mp},
&&\Gamma_{rqm}=\sum^3_{p=0}g_{mp}\,\Gamma^p_{rq}.
\qquad
\mytag{5.11}
\endxalignat
$$
The first formula \mythetag{5.11} is derived from \mythetag{2.4}, \mythetag{2.5}, 
and \mythetag{2.8}. The second formula \mythetag{5.11} is a notation defining
$\Gamma_{rqm}$. In order to transform the second term in the right hand side of 
\mythetag{5.10} we use the formula
$$
\hskip -2em
\sum^2_{\bar s=1}\sum^3_{q=0}L_{\boldsymbol
\Upsilon_{\!r}}\!(G^{k\bar s}_q)\,G^{\,q}_{i\kern 0.5pt\bar s}
=-\sum^2_{\bar s=1}\sum^3_{q=0}G^{k\bar s}_q\,
L_{\boldsymbol\Upsilon_{\!r}}\!(G^{\,q}_{i\kern 0.5pt\bar s})
\mytag{5.12}
$$
This formula is derived by applying the differential operator 
\mythetag{2.15} to the second identity \mythetag{1.1}. Applying 
\mythetag{5.11} and  \mythetag{5.12} to  \mythetag{5.10}, we get
$$
\hskip -2em
\gathered
\sum^2_{k=1}\Alpha^k_{ri}\,d_{kj}=\frac{1}{4}\sum^2_{\bar s=1}
\sum^3_{q=0}\sum^3_{m=0}\sum^2_{\bar r=1}G^{\,m}_{\!j\kern 0.5pt
\bar r}\,G^{\,q}_{i\kern 0.5pt\bar s}\,\Gamma_{rqm}
\,\bd^{\kern 0.5pt\bar r\bar s}\,+\\
+\,\frac{1}{4}\sum^2_{k=1}\sum^2_{\bar s=1}\sum^3_{q=0}
G^{k\bar s}_q\,L_{\boldsymbol\Upsilon_{\!r}}\!(G^{\,q}_{i\kern 0.5pt
\bar s})\,d_{kj}
-\frac{1}{4}\sum^2_{\bar i=1}\sum^2_{\bar j=1}L_{\boldsymbol\Upsilon_{\!r}}
\!(\bd_{\bar j\kern 0.5pt\bar i})\,\bd^{\kern 0.5pt\bar i\bar j}\,d_{ij}.
\endgathered
\mytag{5.13}
$$\par
     Now we replace $m$ by $p$ in the first term and replace $q$ by $p$ 
in the second term in the right hand side of the formula \mythetag{5.13}.
This yields
$$
\hskip -2em
\gathered
\sum^2_{k=1}\Alpha^k_{ri}\,d_{kj}=\frac{1}{4}\sum^2_{\bar s=1}
\sum^3_{q=0}\sum^3_{p=0}\sum^2_{\bar r=1}G^{\,p}_{\!j\kern 0.5pt
\bar r}\,G^{\,q}_{i\kern 0.5pt\bar s}\,\Gamma_{rq\kern 0.3pt p}
\,\bd^{\kern 0.5pt\bar r\bar s}\,+\\
+\,\frac{1}{4}\sum^2_{k=1}\sum^2_{\bar s=1}\sum^3_{p=0}
G^{k\bar s}_p\,L_{\boldsymbol\Upsilon_{\!r}}\!(G^{\,p}_{i\kern 0.5pt
\bar s})\,d_{kj}
-\frac{1}{4}\sum^2_{\bar i=1}\sum^2_{\bar j=1}L_{\boldsymbol\Upsilon_{\!r}}
\!(\bd_{\bar j\kern 0.5pt\bar i})\,\bd^{\kern 0.5pt\bar i\bar j}\,d_{ij}.
\endgathered
\mytag{5.14}
$$
As a result we can apply the first formula \mythetag{5.11} to the second term 
in the right hand side of the formula \mythetag{5.14}. Then we have
$$
\gathered
\sum^2_{k=1}\Alpha^k_{ri}\,d_{kj}=\frac{1}{4}\sum^2_{\bar s=1}
\sum^3_{q=0}\sum^3_{p=0}\sum^2_{\bar r=1}G^{\,p}_{\!j\kern 0.5pt
\bar r}\,G^{\,q}_{i\kern 0.5pt\bar s}\,\Gamma_{rq\kern 0.3pt p}
\,\bd^{\kern 0.5pt\bar r\bar s}\,+\\
+\,\frac{1}{4}\sum^2_{\bar s=1}\sum^2_{\bar r=1}\sum^3_{m=0}\sum^3_{p=0}
G^{\,m}_{\!j\kern 0.5pt\bar r}\,\bd^{\kern 0.5pt\bar r\bar s}\,g_{mp}
\ L_{\boldsymbol\Upsilon_{\!r}}\!(G^{\,p}_{i\kern 0.5pt
\bar s})-\frac{1}{4}\sum^2_{\bar i=1}\sum^2_{\bar j=1}L_{\boldsymbol
\Upsilon_{\!r}}\!(\bd_{\bar j\kern 0.5pt\bar i})\,\bd^{\kern 0.5pt\bar i
\bar j}\,d_{ij}.
\endgathered
\quad
\mytag{5.15}
$$
For the sake of beauty we replace $m$ by $q$ in the formula \mythetag{5.15}:
$$
\gathered
\sum^2_{k=1}\Alpha^k_{ri}\,d_{kj}=\frac{1}{4}\sum^2_{\bar s=1}
\sum^3_{q=0}\sum^3_{p=0}\sum^2_{\bar r=1}G^{\,p}_{\!j\kern 0.5pt
\bar r}\,G^{\,q}_{i\kern 0.5pt\bar s}\,\Gamma_{rq\kern 0.3pt p}
\,\bd^{\kern 0.5pt\bar r\bar s}\,+\\
+\,\frac{1}{4}\sum^2_{\bar s=1}\sum^2_{\bar r=1}\sum^3_{p=0}\sum^3_{q=0}
G^{\,q}_{\!j\kern 0.5pt\bar r}\,\bd^{\kern 0.5pt\bar r\bar s}
\,g_{p\kern 0.5pt q}
\,L_{\boldsymbol\Upsilon_{\!r}}\!(G^{\,p}_{i\kern 0.5pt
\bar s})-\frac{1}{4}\sum^2_{\bar i=1}\sum^2_{\bar j=1}L_{\boldsymbol
\Upsilon_{\!r}}\!(\bd_{\bar j\kern 0.5pt\bar i})\,\bd^{\kern 0.5pt\bar i
\bar j}\,d_{ij}.
\endgathered
\quad
\mytag{5.16}
$$\par
     Now we proceed to the third term in right hand side of the formula 
\mythetag{5.8}. In order to transform it we use the skew symmetry
$d_{i\kern 0.5pt k}=-d_{k\kern 0.5pt i}$:
$$
\hskip -2em
\sum^2_{k=1}\Alpha^k_{rj}\,d_{i\kern 0.5pt k}=
-\sum^2_{k=1}\Alpha^k_{rj}\,d_{k\kern 0.5pt i}.
\mytag{5.17}
$$
Then we can apply the formula \mythetag{5.16} to \mythetag{5.17}. As a result 
we get
$$
\gathered
\sum^2_{k=1}\Alpha^k_{rj}\,d_{i\kern 0.5pt k}
=-\frac{1}{4}\sum^2_{\bar s=1}
\sum^3_{q=0}\sum^3_{p=0}\sum^2_{\bar r=1}G^{\,p}_{\!i\kern 0.5pt
\bar r}\,G^{\,q}_{j\kern 0.5pt\bar s}\,\Gamma_{rq\kern 0.3pt p}
\,\bd^{\kern 0.5pt\bar r\bar s}\,-\\
-\,\frac{1}{4}\sum^2_{\bar s=1}\sum^2_{\bar r=1}\sum^3_{p=0}
\sum^3_{q=0}G^{\,q}_{\!i\kern 0.5pt\bar r}\,\bd^{\kern 0.5pt\bar r
\bar s}\,g_{p\kern 0.5pt q}
\,L_{\boldsymbol\Upsilon_{\!r}}\!(G^{\,p}_{j\kern 0.5pt
\bar s})+\frac{1}{4}\sum^2_{\bar i=1}\sum^2_{\bar j=1}L_{\boldsymbol
\Upsilon_{\!r}}\!(\bd_{\bar j\kern 0.5pt\bar i})\,\bd^{\kern 0.5pt\bar i
\bar j}\,d_{j\kern 0.5pt i}.
\endgathered\quad
\mytag{5.18}
$$
Let's exchange $p$ with $q$ and exchange $\bar r$ with $\bar s$ in
the first term in right hand side of the above formula \mythetag{5.18}.
Moreover, let's do the same in the second term in right hand side of this 
formula and take into account the symmetry $g_{qp}=g_{p\kern 0.7pt q}$ 
and the skew symmetry $\bd^{\kern 0.5pt\bar s\bar r}=-\bd^{\kern 0.5pt
\bar r\bar s}$. As a result of these transformations we obtain 
$$
\gathered
\sum^2_{k=1}\Alpha^k_{rj}\,d_{i\kern 0.5pt k}
=\frac{1}{4}\sum^2_{\bar s=1}
\sum^3_{q=0}\sum^3_{p=0}\sum^2_{\bar r=1}G^{\,q}_{\!i\kern 0.5pt
\bar s}\,G^{\,p}_{j\kern 0.5pt\bar r}\,\Gamma_{rp\kern 0.7pt q}
\,\bd^{\kern 0.5pt\bar r\bar s}\,+\\
+\,\frac{1}{4}\sum^2_{\bar s=1}\sum^2_{\bar r=1}\sum^3_{p=0}
\sum^3_{q=0}L_{\boldsymbol\Upsilon_{\!r}}\!(G^{\,q}_{j\kern 0.5pt
\bar r})\,g_{p\kern 0.7pt q}
\,G^{\,p}_{\!i\kern 0.5pt\bar s}\,\bd^{\kern 0.5pt\bar r\bar s}
-\frac{1}{4}\sum^2_{\bar i=1}\sum^2_{\bar j=1}L_{\boldsymbol
\Upsilon_{\!r}}\!(\bd_{\bar j\kern 0.5pt\bar i})\,\bd^{\kern 0.5pt\bar i
\bar j}\,d_{ij}.
\endgathered\quad
\mytag{5.19}
$$
Now let's add the equalities \mythetag{5.16} and  \mythetag{5.19}. This
yields
$$
\gathered
\sum^2_{k=1}\Alpha^k_{ri}\,d_{kj}+\sum^2_{k=1}\Alpha^k_{rj}
\,d_{i\kern 0.5pt k}
=\frac{1}{4}\sum^2_{\bar s=1}
\sum^3_{q=0}\sum^3_{p=0}\sum^2_{\bar r=1}G^{\,p}_{\!j\kern 0.5pt
\bar r}\,G^{\,q}_{i\kern 0.5pt\bar s}\bigl(\Gamma_{rq\kern 0.3pt p}
\,+\\
+\,\Gamma_{rp\kern 0.7pt q}\bigr)\,\bd^{\kern 0.5pt\bar r\bar s}
+\frac{1}{4}\sum^2_{\bar s=1}\sum^2_{\bar r=1}\sum^3_{p=0}\sum^3_{q=0}
\bigl(L_{\boldsymbol\Upsilon_{\!r}}\!(G^{\,q}_{j\kern 0.5pt
\bar r})\,G^{\,p}_{\!i\kern 0.5pt\bar s}+G^{\,q}_{\!j\kern 0.5pt
\bar r}\,L_{\boldsymbol\Upsilon_{\!r}}\!(G^{\,p}_{i\kern 0.5pt\bar s})
\bigr)\,\times\\
\times\,\bd^{\kern 0.5pt\bar r\bar s}\,g_{p\kern 0.5pt q}
-\frac{1}{2}\sum^2_{\bar i=1}\sum^2_{\bar j=1}L_{\boldsymbol
\Upsilon_{\!r}}\!(\bd_{\bar j\kern 0.5pt\bar i})\,\bd^{\kern 0.5pt\bar i
\bar j}\,d_{ij}.
\endgathered
\quad
\mytag{5.20}
$$
Let's recall that the differential operator \mythetag{2.15}, when applied 
to a product, obeys the Leibniz rule. Then we can perform the following 
transformations:
$$
\gather
\sum^2_{\bar s=1}\sum^2_{\bar r=1}\sum^3_{p=0}\sum^3_{q=0}
\bigl(L_{\boldsymbol\Upsilon_{\!r}}\!(G^{\,q}_{j\kern 0.5pt
\bar r})\,G^{\,p}_{\!i\kern 0.5pt\bar s}+G^{\,q}_{\!j\kern 0.5pt
\bar r}\,L_{\boldsymbol\Upsilon_{\!r}}\!(G^{\,p}_{i\kern 0.5pt\bar s})
\bigr)\,\bd^{\kern 0.5pt\bar r\bar s}\,g_{p\kern 0.5pt q}=\\
=\sum^2_{\bar s=1}\sum^2_{\bar r=1}\sum^3_{p=0}\sum^3_{q=0}
L_{\boldsymbol\Upsilon_{\!r}}\!(G^{\,q}_{j\kern 0.5pt
\bar r}\,G^{\,p}_{\!i\kern 0.5pt\bar s})\,\bd^{\kern 0.5pt\bar r\bar s}
\,g_{p\kern 0.5pt q}=\sum^2_{\bar s=1}\sum^2_{\bar r=1}\sum^3_{p=0}
\sum^3_{q=0}L_{\boldsymbol\Upsilon_{\!r}}\!(G^{\,q}_{j\kern 0.5pt
\bar r}\,\times\\
\times\,G^{\,p}_{\!i\kern 0.5pt\bar s}\ g_{p\kern 0.5pt q})
\,\bd^{\kern 0.5pt\bar r\bar s}-\sum^2_{\bar s=1}\sum^2_{\bar r=1}
\sum^3_{p=0}\sum^3_{q=0}G^{\,q}_{j\kern 0.5pt\bar r}
\,G^{\,p}_{\!i\kern 0.5pt\bar s}
\,L_{\boldsymbol\Upsilon_{\!r}}\!(g_{p\kern 0.5pt q})
\,\bd^{\kern 0.5pt\bar r\bar s}.
\endgather
$$
Due to \mythetag{2.3} and \mythetag{2.8} the first formula \mythetag{5.11}
yields
$$
\hskip -2em
\sum^2_{\bar a=1}\sum^2_{k=1}G^{k\bar a}_p\,d_{kj}\,\bd_{\bar a\kern 0.3pt
\bar r}=\sum^3_{q=0}G^{\,q}_{\!j\kern 0.5pt\bar r}\,g_{p\kern 0.5pt q}.
\mytag{5.21}
$$
Applying the relationship \mythetag{5.21}, we can continue the above 
transformations
$$
\gather
\sum^2_{\bar s=1}\sum^2_{\bar r=1}\sum^3_{p=0}\sum^3_{q=0}
\bigl(L_{\boldsymbol\Upsilon_{\!r}}\!(G^{\,q}_{j\kern 0.5pt
\bar r})\,G^{\,p}_{\!i\kern 0.5pt\bar s}+G^{\,q}_{\!j\kern 0.5pt
\bar r}\,L_{\boldsymbol\Upsilon_{\!r}}\!(G^{\,p}_{i\kern 0.5pt\bar s})
\bigr)\,\bd^{\kern 0.5pt\bar r\bar s}\,g_{p\kern 0.5pt q}=\\
=\sum^2_{k=1}\sum^2_{\bar a=1}\sum^2_{\bar s=1}\sum^2_{\bar r=1}
\sum^3_{p=0}L_{\boldsymbol\Upsilon_{\!r}}\!(
G^{k\bar a}_p\,G^{\,p}_{\!i\kern 0.5pt\bar s}
\,d_{kj}\,\bd_{\bar a\kern 0.3pt\bar r})
\,\bd^{\kern 0.5pt\bar r\bar s}\,-\kern 4em\\
\kern 8em -\,\sum^2_{\bar s=1}\sum^2_{\bar r=1}
\sum^3_{p=0}\sum^3_{q=0}G^{\,q}_{j\kern 0.5pt\bar r}
\,G^{\,p}_{\!i\kern 0.5pt\bar s}
\,L_{\boldsymbol\Upsilon_{\!r}}\!(g_{p\kern 0.5pt q})
\,\bd^{\kern 0.5pt\bar r\bar s}.
\endgather
$$
Now we can apply the second identity \mythetag{1.1} to the first term
in right hand side of the above equality. As a result we derive 
$$
\gathered
\sum^2_{\bar s=1}\sum^2_{\bar r=1}\sum^3_{p=0}\sum^3_{q=0}
\bigl(L_{\boldsymbol\Upsilon_{\!r}}\!(G^{\,q}_{j\kern 0.5pt
\bar r})\,G^{\,p}_{\!i\kern 0.5pt\bar s}+G^{\,q}_{\!j\kern 0.5pt
\bar r}\,L_{\boldsymbol\Upsilon_{\!r}}\!(G^{\,p}_{i\kern 0.5pt\bar s})
\bigr)\,\bd^{\kern 0.5pt\bar r\bar s}\,g_{p\kern 0.5pt q}=\\
=\sum^2_{\bar s=1}\sum^2_{\bar r=1}L_{\boldsymbol\Upsilon_{\!r}}\!(
2\,d_{ij}\,\bd_{\bar s\kern 0.3pt\bar r})\,\bd^{\kern 0.5pt\bar r\bar s}
-\sum^2_{\bar s=1}\sum^2_{\bar r=1}
\sum^3_{p=0}\sum^3_{q=0}G^{\,q}_{j\kern 0.5pt\bar r}
\,G^{\,p}_{\!i\kern 0.5pt\bar s}
\,L_{\boldsymbol\Upsilon_{\!r}}\!(g_{p\kern 0.5pt q})
\,\bd^{\kern 0.5pt\bar r\bar s}.
\endgathered
\mytag{5.22}
$$
Applying the Leibniz rule to the first term in the right hand side 
of \mythetag{5.22} and then substituting \mythetag{5.22} back into 
\mythetag{5.20}, we obtain
$$
\gathered
\sum^2_{k=1}\Alpha^k_{ri}\,d_{kj}+\sum^2_{k=1}\Alpha^k_{rj}
\,d_{i\kern 0.5pt k}
=\frac{1}{4}\sum^2_{\bar s=1}
\sum^3_{q=0}\sum^3_{p=0}\sum^2_{\bar r=1}G^{\,p}_{\!j\kern 0.5pt
\bar r}\,G^{\,q}_{i\kern 0.5pt\bar s}\bigl(\Gamma_{rq\kern 0.3pt p}
\,+\\
+\,\Gamma_{rp\kern 0.7pt q}\bigr)\,\bd^{\kern 0.5pt\bar r\bar s}
+L_{\boldsymbol\Upsilon_{\!r}}\!(d_{ij})
-\frac{1}{4}
\sum^2_{\bar s=1}\sum^2_{\bar r=1}
\sum^3_{p=0}\sum^3_{q=0}G^{\,q}_{j\kern 0.5pt\bar r}
\,G^{\,p}_{\!i\kern 0.5pt\bar s}
\,L_{\boldsymbol\Upsilon_{\!r}}\!(g_{p\kern 0.5pt q})
\,\bd^{\kern 0.5pt\bar r\bar s}.
\endgathered
\quad
\mytag{5.23}
$$\par
     The next step is to transform $\Gamma_{rq\kern 0.3pt p}$ in 
\mythetag{5.23}. For this purpose we substitute \mythetag{5.4} into
the second formula \mythetag{5.11}. As a result we get
$$
\hskip -2em
\gathered
\Gamma_{rq\kern 0.3pt p}=\frac{1}{2}
\left(L_{\boldsymbol\Upsilon_{\!r}}\!(g_{p\kern 0.5pt q})
+L_{\boldsymbol\Upsilon_{\!q}}\!(g_{rp})
-L_{\boldsymbol\Upsilon_{\!p}}\!(g_{rq})\right)+\\
+\sum^3_{s=0}\frac{c^{\,s}_{rq}}{2}\,g_{sp}
-\sum^3_{s=0}\frac{c^{\,s}_{rp}}{2}
\,g_{sq}-\sum^3_{s=0}\frac{c^{\,s}_{q\kern 0.5pt p}}{2}
\,g_{s\kern 0.5pt r}.
\endgathered
\mytag{5.24}
$$
When symmetrizing with respect to the indices $p$ and $q$ most of the 
terms in \mythetag{5.24} do cancel each other. The rest of them yield
$$
\hskip -2em
\Gamma_{rq\kern 0.3pt p}+\Gamma_{rp\kern 0.7pt q}
=L_{\boldsymbol\Upsilon_{\!r}}\!(g_{p\kern 0.5pt q}).
\mytag{5.25}
$$
Substituting \mythetag{5.25} back into the formula \mythetag{5.23},
we derive
$$
\sum^2_{k=1}\Alpha^k_{ri}\,d_{kj}+\sum^2_{k=1}\Alpha^k_{rj}
\,d_{i\kern 0.5pt k}=L_{\boldsymbol\Upsilon_{\!r}}\!(d_{ij}).
\mytag{5.26}
$$
Comparing \mythetag{5.26} with \mythetag{5.8}, we find that these
formulas are equivalent. Thus we conclude that for the spinor connection
with the components \mythetag{1.3} the first concordance condition 
\mythetag{5.7} is fulfilled.\par
     Now let's proceed to verifying the second concordance condition 
\mythetag{5.7}. Its coordinate presentation is given by the formula
\mythetag{5.9}. Let's begin with the second term in \mythetag{5.9}.
Applying \mythetag{1.3} to this term, we derive
$$
\gathered
\sum^2_{k=1}\Alpha^i_{\kern 0.5pt rk}\,G^{k\kern 0.5pt \bar i}_p
=\sum^2_{k=1}\sum^2_{\bar s=1}\sum^3_{q=0}\sum^3_{m=0}
\frac{G^{i\bar s}_q\,\Gamma^q_{rm}
\,G^{\,m}_{k\kern 0.5pt\bar s}\,G^{k\kern 0.5pt \bar i}_p}{4}\,-\\
-\sum^2_{k=1}\sum^2_{\bar s=1}\sum^3_{m=0}
\frac{L_{\boldsymbol\Upsilon_{\!r}}\!(G^{i\bar s}_m)
\,G^{\,m}_{k\kern 0.5pt\bar s}\,G^{k\kern 0.5pt \bar i}_p}{4}
-\sum^2_{\bar u=1}\sum^2_{\bar s=1}
\frac{L_{\boldsymbol\Upsilon_{\!r}}\!(\bd_{\bar s\kern 0.5pt\bar u})
\,\bd^{\kern 0.5pt\bar u\kern 0.3pt\bar s}
\,G^{i\kern 0.5pt \bar i}_p}{4}.
\endgathered
\qquad
\mytag{5.27}
$$
In a similar way, applying \mythetag{1.3} to the third term in
\mythetag{5.9}, we obtain
$$
\gathered
\sum^2_{\bar k=1}\overline{\Alpha^{\bar i}_{r\bar k}}
\,G^{i\kern 0.5pt \bar k}_p
=\sum^2_{\bar k=1}\sum^2_{s=1}\sum^3_{q=0}\sum^3_{m=0}
\frac{\overline{G^{\raise 0.7pt\hbox{$\ssize\bar is$}}_q}
\,\Gamma^q_{rm}
\,\overline{G\raise 5	pt\hbox{$\kern 1pt\ssize m$}
\lower 2.6pt\hbox{$\kern -8pt\ssize \bar k\kern 0.5pt s$}}
\,G^{i\kern 0.5pt \bar k}_p}{4}\,-\\
-\sum^2_{\bar k=1}\sum^2_{s=1}\sum^3_{m=0}
\frac{L_{\boldsymbol\Upsilon_{\!r}}\!(\overline{G^{\bar is}_m})
\,\overline{G\raise 5	pt\hbox{$\kern 1pt\ssize m$}
\lower 2.6pt\hbox{$\kern -8pt\ssize \bar k\kern 0.5pt s$}}
\,G^{i\kern 0.5pt \bar k}_p}{4}
-\sum^2_{u=1}\sum^2_{s=1}
\frac{L_{\boldsymbol\Upsilon_{\!r}}\!
(\overline{\bd_{s\kern 0.3pt u}})
\,\overline{\bd^{\kern 0.5pt us}}
\,G^{i\kern 0.5pt \bar i}_p}{4}.
\endgathered
\qquad
\mytag{5.28}
$$
In order to transform \mythetag{5.28} we apply \mythetag{2.6},
\mythetag{2.7}, and \mythetag{2.9}. This yields
$$
\gathered
\sum^2_{\bar k=1}\overline{\Alpha^{\bar i}_{r\bar k}}
\,G^{i\kern 0.5pt \bar k}_p
=\sum^2_{\bar k=1}\sum^2_{s=1}\sum^3_{q=0}\sum^3_{m=0}
\frac{G^{s\bar i}_q\,\Gamma^q_{rm}
\,G\raise 5	pt\hbox{$\kern 1pt\ssize m$}
\lower 2.6pt\hbox{$\kern -8pt\ssize s\bar k$}
\ G^{i\kern 0.5pt \bar k}_p}{4}\,-\\
-\sum^2_{\bar k=1}\sum^2_{s=1}\sum^3_{m=0}
\frac{L_{\boldsymbol\Upsilon_{\!r}}\!(G^{s\bar i}_m)
\,G\raise 5	pt\hbox{$\kern 1pt\ssize m$}
\lower 2.6pt\hbox{$\kern -9pt\ssize s\bar k$}
\ G^{i\kern 0.5pt \bar k}_p}{4}
-\sum^2_{u=1}\sum^2_{s=1}
\frac{L_{\boldsymbol\Upsilon_{\!r}}\!
(d_{s\kern 0.3pt u})
\,d^{\kern 0.5pt us}
\,G^{i\kern 0.5pt \bar i}_p}{4}.
\endgathered
\qquad
\mytag{5.29}
$$
Now we replace $s$ by $k$ and replace $\bar k$ by $\bar s$ in 
the right hand side of \mythetag{5.29}:
$$
\gathered
\sum^2_{\bar k=1}\overline{\Alpha^{\bar i}_{r\bar k}}
\,G^{i\kern 0.5pt \bar k}_p
=\sum^2_{k=1}\sum^2_{\bar s=1}\sum^3_{q=0}\sum^3_{m=0}
\frac{G^{i\kern 0.5pt \bar s}_p\,\Gamma^q_{rm}\,G^{\,m}_{k\bar s}
\,G^{k\bar i}_q}{4}\,-\\
-\sum^2_{k=1}\sum^2_{\bar s=1}\sum^3_{m=0}
\frac{G^{i\kern 0.5pt \bar s}_p\,G^{\,m}_{k\bar s}
\,L_{\boldsymbol\Upsilon_{\!r}}\!(G^{k\bar i}_m)}{4}
-\sum^2_{u=1}\sum^2_{s=1}
\frac{L_{\boldsymbol\Upsilon_{\!r}}\!
(d_{s\kern 0.3pt u})\,d^{\kern 0.5pt us}
\,G^{i\kern 0.5pt \bar i}_p}{4}.
\endgathered
\qquad
\mytag{5.30}
$$
Let's return back to the formula \mythetag{5.12} and rewrite it
as follows:
$$
\hskip -2em
\sum^2_{\bar s=1}\sum^3_{m=0}L_{\boldsymbol
\Upsilon_{\!r}}\!(G^{i\bar s}_m)\,G^{\,m}_{k\kern 0.5pt\bar s}
=-\sum^2_{\bar s=1}\sum^3_{m=0}G^{i\bar s}_m\,
L_{\boldsymbol\Upsilon_{\!r}}\!(G^{\,m}_{k\kern 0.5pt\bar s}).
\mytag{5.31}
$$
There is a formula very similar to \mythetag{5.31}. Here it is:
$$
\hskip -2em
\sum^2_{k=1}\sum^3_{m=0}G^{\,m}_{k\bar s}
\,L_{\boldsymbol\Upsilon_{\!r}}\!(G^{k\bar i}_m)
=-\sum^2_{k=1}\sum^3_{m=0}
L_{\boldsymbol\Upsilon_{\!r}}\!(G^{\,m}_{k\bar s})
\,G^{k\bar i}_m.
\mytag{5.32}
$$
We apply \mythetag{5.31} to \mythetag{5.27} and apply \mythetag{5.32} 
to \mythetag{5.30}. Then we add these two formulas. As a result we
obtain the following formula:
$$
\gathered
\sum^2_{k=1}\Alpha^i_{\kern 0.5pt rk}\,G^{k\kern 0.5pt \bar i}_p
+\sum^2_{\bar k=1}\overline{\Alpha^{\bar i}_{r\bar k}}
\,G^{i\kern 0.5pt \bar k}_p
=\sum^2_{k=1}\sum^2_{\bar s=1}\sum^3_{q=0}\sum^3_{m=0}
\frac{\bigl(G^{i\bar s}_q\,G^{k\kern 0.5pt \bar i}_p
+G^{i\kern 0.5pt \bar s}_p\,\,G^{k\bar i}_q\bigr)}{4}
\,\times\\
\times\,\Gamma^q_{rm}\,G^{\,m}_{k\kern 0.5pt\bar s}
+\sum^2_{k=1}\sum^2_{\bar s=1}\sum^3_{m=0}
\frac{\bigl(G^{i\bar s}_m\,G^{k\kern 0.5pt \bar i}_p
+G^{i\kern 0.5pt \bar s}_p\,G^{k\bar i}_m\bigr)}{4}
\,L_{\boldsymbol\Upsilon_{\!r}}\!(G^{\,m}_{k\kern 0.5pt\bar s})\,-\\
-\sum^2_{\bar u=1}\sum^2_{\bar s=1}
\frac{L_{\boldsymbol\Upsilon_{\!r}}\!(\bd_{\bar s\kern 0.3pt\bar u})
\,\bd^{\kern 0.5pt\bar u\kern 0.3pt\bar s}
\,G^{i\kern 0.5pt \bar i}_p}{4}
-\sum^2_{u=1}\sum^2_{s=1}
\frac{L_{\boldsymbol\Upsilon_{\!r}}\!(d_{s\kern 0.3pt u})
\,d^{\kern 0.5pt us}\,G^{i\kern 0.5pt \bar i}_p}{4}.
\endgathered
\quad
\mytag{5.33}
$$
In order to transform \mythetag{5.33} we use the formula \mythetag{4.2}
derived from the identity \mythetag{1.2}. We rewrite this formula as
follows:
$$
\hskip -2em
\aligned
&G^{i\bar s}_q\,G^{k\kern 0.5pt \bar i}_p
+G^{i\kern 0.5pt \bar s}_p\,G^{k\bar i}_q
=G^{i\kern 0.5pt\bar i}_q\,G^{k\bar s}_p
+G^{i\kern 0.5pt \bar i}_p\,G^{k\bar s}_q
-2\,d^{\kern 0.5pt ik}\,\bd^{\kern 0.5pt \bar i\bar s}
\,g_{p\kern 0.7pt q},\\
\vspace{1ex}
&G^{i\bar s}_m\,G^{k\kern 0.5pt \bar i}_p
+G^{i\kern 0.5pt \bar s}_p\,G^{k\bar i}_m
=G^{i\kern 0.5pt\bar i}_m\,G^{k\bar s}_p
+G^{i\kern 0.5pt \bar i}_p\,G^{k\bar s}_m
-2\,d^{\kern 0.5pt ik}\,\bd^{\kern 0.5pt \bar i\bar s}
\,g_{p\kern 0.7pt m}.
\endaligned
\mytag{5.34}
$$
Substituting \mythetag{5.34} back into the formula \mythetag{5.33}, 
we derive
$$
\gathered
\sum^2_{k=1}\Alpha^i_{\kern 0.5pt rk}\,G^{k\kern 0.5pt \bar i}_p
+\sum^2_{\bar k=1}\overline{\Alpha^{\bar i}_{r\bar k}}
\,G^{i\kern 0.5pt \bar k}_p
=\sum^2_{k=1}\sum^2_{\bar s=1}\sum^3_{q=0}\sum^3_{m=0}
\frac{\bigl(G^{i\kern 0.5pt\bar i}_q\,G^{k\bar s}_p
+G^{i\kern 0.5pt \bar i}_p\,G^{k\bar s}_q\bigr)}{4}
\,\times\\
\times\,\Gamma^q_{rm}\,G^{\,m}_{k\kern 0.5pt\bar s}
+\sum^2_{k=1}\sum^2_{\bar s=1}\sum^3_{m=0}
\frac{\bigl(G^{i\kern 0.5pt\bar i}_m\,G^{k\bar s}_p
+G^{i\kern 0.5pt \bar i}_p\,G^{k\bar s}_m\bigr)}{4}
\,L_{\boldsymbol\Upsilon_{\!r}}\!(G^{\,m}_{k\kern 0.5pt\bar s})\,-\\
-\sum^2_{k=1}\sum^2_{\bar s=1}\sum^3_{q=0}\sum^3_{m=0}
\frac{d^{\kern 0.5pt ik}\,\bd^{\kern 0.5pt \bar i\bar s}
\,g_{p\kern 0.7pt q}}{2}
\,\Gamma^q_{rm}\,G^{\,m}_{k\kern 0.5pt\bar s}
-\sum^2_{k=1}\sum^2_{\bar s=1}\sum^3_{m=0}
\frac{d^{\kern 0.5pt ik}\,\bd^{\kern 0.5pt \bar i\bar s}
\,g_{p\kern 0.7pt m}}{2}\,L_{\boldsymbol\Upsilon_{\!r}}\!(
G^{\,m}_{k\kern 0.5pt\bar s})\,-\\
-\sum^2_{\bar u=1}\sum^2_{\bar s=1}
\frac{L_{\boldsymbol\Upsilon_{\!r}}\!(\bd_{\bar s\kern 0.3pt\bar u})
\,\bd^{\kern 0.5pt\bar u\kern 0.3pt\bar s}
\,G^{i\kern 0.5pt \bar i}_p}{4}
-\sum^2_{u=1}\sum^2_{s=1}
\frac{L_{\boldsymbol\Upsilon_{\!r}}\!(d_{s\kern 0.3pt u})
\,d^{\kern 0.5pt us}\,G^{i\kern 0.5pt \bar i}_p}{4}.
\endgathered
$$
In order to transform this equality we use the first identity 
\mythetag{1.1} and the second formula \mythetag{5.11}. As a result
we obtain the equality
$$
\gathered
\sum^2_{k=1}\Alpha^i_{\kern 0.5pt rk}\,G^{k\kern 0.5pt \bar i}_p
+\sum^2_{\bar k=1}\overline{\Alpha^{\bar i}_{r\bar k}}
\,G^{i\kern 0.5pt \bar k}_p
=\sum^3_{q=0}\frac{\bigl(G^{i\kern 0.5pt\bar i}_q\,\Gamma^q_{rp}
+G^{i\kern 0.5pt \bar i}_p\,\,\Gamma^q_{rq}\bigr)}{2}\,+\\
+\sum^2_{k=1}\sum^2_{\bar s=1}\sum^3_{m=0}
\frac{G^{i\kern 0.5pt\bar i}_m\,G^{k\bar s}_p}{4}
\,L_{\boldsymbol\Upsilon_{\!r}}\!(G^{\,m}_{k\kern 0.5pt\bar s})
+\sum^2_{k=1}\sum^2_{\bar s=1}\sum^3_{m=0}
\frac{G^{i\kern 0.5pt \bar i}_p\,G^{k\bar s}_m}{4}
\,L_{\boldsymbol\Upsilon_{\!r}}\!(G^{\,m}_{k\kern 0.5pt\bar s})\,-\\
-\sum^2_{k=1}\sum^2_{\bar s=1}\sum^3_{m=0}
\frac{d^{\kern 0.5pt ik}\,\bd^{\kern 0.5pt \bar i\bar s}}{2}
\,\Gamma_{rmp}\,G^{\,m}_{k\kern 0.5pt\bar s}
-\sum^2_{k=1}\sum^2_{\bar s=1}\sum^3_{m=0}
\frac{d^{\kern 0.5pt ik}\,\bd^{\kern 0.5pt \bar i\bar s}
\,g_{p\kern 0.7pt m}}{2}\,L_{\boldsymbol\Upsilon_{\!r}}\!(
G^{\,m}_{k\kern 0.5pt\bar s})\,-\\
-\sum^2_{\bar u=1}\sum^2_{\bar s=1}
\frac{L_{\boldsymbol\Upsilon_{\!r}}\!(\bd_{\bar s\kern 0.3pt\bar u})
\,\bd^{\kern 0.5pt\bar u\kern 0.3pt\bar s}
\,G^{i\kern 0.5pt \bar i}_p}{4}
-\sum^2_{u=1}\sum^2_{s=1}
\frac{L_{\boldsymbol\Upsilon_{\!r}}\!(d_{s\kern 0.3pt u})
\,d^{\kern 0.5pt us}\,G^{i\kern 0.5pt \bar i}_p}{4}.
\endgathered
\qquad
\mytag{5.35}
$$
Applying the differential operator \mythetag{2.15} to the second identity
\mythetag{1.1}, we derive
$$
\hskip -2em
\sum^3_{m=0}G^{i\kern 0.5pt\bar i}_m\,L_{\boldsymbol\Upsilon_{\!r}}\!(
G^{\,m}_{k\kern 0.5pt\bar s})=-\sum^3_{m=0}L_{\boldsymbol\Upsilon_{\!r}}\!(
G^{i\kern 0.5pt\bar i}_m)\,G^{\,m}_{k\kern 0.5pt\bar s}.
\mytag{5.36}
$$
The equality \mythetag{5.36} is similar to \mythetag{5.31} and \mythetag{5.32}.
Now we use the second identity \mythetag{1.1} once more in order to perform
the following calculations:
$$
\gather
0=\sum^2_{k=1}\sum^2_{\bar s=1}\sum^3_{m=0}
L_{\boldsymbol\Upsilon_{\!r}}\!(G^{k\bar s}_m
\,G^{\,m}_{k\kern 0.5pt\bar s})=
\sum^2_{k=1}\sum^2_{\bar s=1}\sum^3_{m=0}
\sum^2_{q=1}\sum^2_{\bar r=1}\sum^3_{n=0}
L_{\boldsymbol\Upsilon_{\!r}}\!(G^{\,n}_{q\kern 0.5pt\bar r}
\,\times\\
\times\,d^{\kern 0.5ptqk}\,\bd^{\kern 0.5pt\bar r\bar s}
g_{nm}\,G^{\,m}_{k\kern 0.5pt\bar s})
=\sum^2_{k=1}\sum^2_{\bar s=1}\sum^3_{m=0}
\sum^2_{q=1}\sum^2_{\bar r=1}\sum^3_{n=0}
L_{\boldsymbol\Upsilon_{\!r}}\!(G^{\,n}_{q\kern 0.5pt\bar r})
\,d^{\kern 0.5ptqk}\,\bd^{\kern 0.5pt\bar r\bar s}\,g_{nm}
\,\times\\
\times\,G^{\,m}_{k\kern 0.5pt\bar s}
+\sum^2_{k=1}\sum^2_{\bar s=1}\sum^3_{m=0}
\sum^2_{q=1}\sum^2_{\bar r=1}\sum^3_{n=0}
G^{\,n}_{q\kern 0.5pt\bar r}\,
L_{\boldsymbol\Upsilon_{\!r}}\!(d^{\kern 0.5ptqk}
\,\bd^{\kern 0.5pt\bar r\bar s}\,g_{nm})
\,G^{\,m}_{k\kern 0.5pt\bar s}\,+\\
+\,\sum^2_{k=1}\sum^2_{\bar s=1}\sum^3_{m=0}
\sum^2_{q=1}\sum^2_{\bar r=1}\sum^3_{n=0}
G^{\,n}_{q\kern 0.5pt\bar r}\,d^{\kern 0.5ptqk}
\,\bd^{\kern 0.5pt\bar r\bar s}\,g_{nm}\,
L_{\boldsymbol\Upsilon_{\!r}}\!(G^{\,m}_{k\kern 0.5pt\bar s})\,=\\
=\sum^2_{k=1}\sum^2_{\bar s=1}\sum^3_{m=0}
\sum^2_{q=1}\sum^2_{\bar r=1}\sum^3_{n=0}
G^{\,n}_{q\kern 0.5pt\bar r}\,
L_{\boldsymbol\Upsilon_{\!r}}\!(d^{\kern 0.5ptqk}
\,\bd^{\kern 0.5pt\bar r\bar s}\,g_{nm})
\,G^{\,m}_{k\kern 0.5pt\bar s}\,+\\
+\sum^2_{q=1}\sum^2_{\bar r=1}\sum^3_{n=0}
L_{\boldsymbol\Upsilon_{\!r}}\!(G^{\,n}_{q\kern 0.5pt\bar r})
\,G^{q\bar r}_n+\sum^2_{k=1}\sum^2_{\bar s=1}\sum^3_{m=0}
G^{k\bar s}_m\,L_{\boldsymbol\Upsilon_{\!r}}\!(G^{\,m}_{k
\kern 0.5pt\bar s}).
\endgather
$$
It easy to see that the last two terms in the above formula are equal
to each other. For this reason the above formula is equivalent to
$$
\sum^2_{k=1}\sum^2_{\bar s=1}\sum^3_{m=0}
G^{k\bar s}_m\,L_{\boldsymbol\Upsilon_{\!r}}\!(G^{\,m}_{k
\kern 0.5pt\bar s})=-\frac{1}{2}\sum^2_{k=1}\sum^2_{\bar s=1}\sum^3_{m=0}
\sum^2_{q=1}\sum^2_{\bar r=1}\sum^3_{n=0}G^{\,n}_{q\kern 0.5pt\bar r}\,
L_{\boldsymbol\Upsilon_{\!r}}\!(d^{\kern 0.5ptqk}
\,\bd^{\kern 0.5pt\bar r\bar s}\,g_{nm})\,G^{\,m}_{k\kern 0.5pt\bar s}.
$$
The operator $L_{\boldsymbol\Upsilon_{\!r}}$ obeys the Leibniz rule.
Therefore, we have
$$
\gathered
\sum^2_{k=1}\sum^2_{\bar s=1}\sum^3_{m=0}
G^{k\bar s}_m\,L_{\boldsymbol\Upsilon_{\!r}}\!(G^{\,m}_{k
\kern 0.5pt\bar s})=
-\frac{1}{2}\sum^2_{k=1}\sum^2_{\bar s=1}\sum^3_{m=0}
\sum^2_{q=1}\sum^2_{\bar r=1}\sum^3_{n=0}G^{\,n}_{q\kern 0.5pt\bar r}
\,L_{\boldsymbol\Upsilon_{\!r}}\!(d^{\kern 0.5ptqk})
\,\times\\
\times\,
\bd^{\kern 0.5pt\bar r\bar s}\,g_{nm}\,G^{\,m}_{k\kern 0.5pt\bar s}
-\frac{1}{2}\sum^2_{k=1}\sum^2_{\bar s=1}\sum^3_{m=0}
\sum^2_{q=1}\sum^2_{\bar r=1}\sum^3_{n=0}G^{\,n}_{q\kern 0.5pt\bar r}\,
d^{\kern 0.5pt qk}\,L_{\boldsymbol\Upsilon_{\!r}}\!(\bd^{\kern 0.5pt
\bar r\bar s})\,\times\\
\times\,g_{nm}\,G^{\,m}_{k\kern 0.5pt\bar s}
-\frac{1}{2}\sum^2_{k=1}\sum^2_{\bar s=1}\sum^3_{m=0}
\sum^2_{q=1}\sum^2_{\bar r=1}\sum^3_{n=0}
G^{\,n}_{q\kern 0.5pt\bar r}\,d^{\kern 0.5pt qk}
\,\bd^{\kern 0.5pt\bar r\bar s}\,L_{\boldsymbol\Upsilon_{\!r}}\!(\,g_{nm})
\,G^{\,m}_{k\kern 0.5pt\bar s}.
\endgathered
\quad
\mytag{5.37}
$$
Note that the first formula \mythetag{5.11}, can be rewritten as follows:
$$
\hskip -2em
\sum^2_{\bar s=1}\sum^3_{m=0}\bd^{\kern 0.5pt\bar r\bar s}\,g_{nm}\,
G^{\,m}_{k\kern 0.5pt\bar s}=\sum^2_{s=1}G^{s\kern 0.5pt\bar r}_n\,d_{ks}.
\mytag{5.38}
$$
The equality \mythetag{5.11} itself was derived from \mythetag{2.3}, 
\mythetag{2.4}, \mythetag{2.5}, and \mythetag{2.8}. By analogy to it
one can derive the following two equalities:
$$
\align
&\hskip -2em
\sum^2_{k=1}\sum^3_{m=0}d^{\kern 0.5pt qk}\,g_{nm}
\,G^{\,m}_{k\kern 0.5pt\bar s}=\sum^2_{\bar k=1}
G^{q\kern 0.5pt\bar k}_n\,\bd_{\bar s\bar k},
\mytag{5.39}\\
&\hskip -2em
\sum^2_{k=1}\sum^2_{\bar s=1}
d^{\kern 0.5pt qk}\bd^{\kern 0.5pt\bar r\bar s}
\,G^{\,m}_{k\kern 0.5pt\bar s}=
\sum^3_{p=0}G^{q\kern 0.5pt\bar r}_p\,g^{p\kern 0.5pt m}.
\mytag{5.40}
\endalign
$$
Now we apply the formulas \mythetag{5.38}, \mythetag{5.39}, and 
\mythetag{5.40} to the formula \mythetag{5.37} and take into account 
the identities \mythetag{1.1}. As a result we obtain
$$
\gathered
\sum^2_{k=1}\sum^2_{\bar s=1}\sum^3_{m=0}
G^{k\bar s}_m\,L_{\boldsymbol\Upsilon_{\!r}}\!(G^{\,m}_{k
\kern 0.5pt\bar s})=-2\sum^2_{k=1}\sum^2_{q=1}
\,L_{\boldsymbol\Upsilon_{\!r}}\!(d^{\kern 0.5ptqk})
\,d_{kq}\,-\\
-\,2\sum^2_{\bar r=1}\sum^2_{\bar s=1}
L_{\boldsymbol\Upsilon_{\!r}}\!(\bd^{\kern 0.5pt
\bar r\bar s})\,\bd_{\bar s\bar r}-\sum^3_{m=0}\sum^3_{n=0}
g^{nm}\,L_{\boldsymbol\Upsilon_{\!r}}\!(\,g_{nm}).
\endgathered
\quad
\mytag{5.41}
$$
From \mythetag{2.4} and \mythetag{2.8} by applying the differential 
operator $L_{\boldsymbol\Upsilon_{\!r}}$ to them we derive
$$
\hskip -2em
\aligned
&\sum^2_{k=1}\sum^2_{q=1}
L_{\boldsymbol\Upsilon_{\!r}}\!(d^{\kern 0.5pt qk})
\,d_{kq}=-\sum^2_{k=1}\sum^2_{q=1}d^{\kern 0.5pt qk}
\,L_{\boldsymbol\Upsilon_{\!r}}\!(d_{kq}),\\
&\sum^2_{\bar r=1}\sum^2_{\bar s=1}
L_{\boldsymbol\Upsilon_{\!r}}\!(\bd^{\kern 0.5pt
\bar r\bar s})\,\bd_{\bar s\bar r}=
-\sum^2_{\bar r=1}\sum^2_{\bar s=1}
\bd^{\kern 0.5pt\bar r\bar s}
\,L_{\boldsymbol\Upsilon_{\!r}}\!(\bd_{\bar s\bar r}).
\endaligned
\mytag{5.42}
$$
Substituting \mythetag{5.42} into \mythetag{5.41}, we transform 
\mythetag{5.41} as follows:
$$
\gathered
\sum^2_{k=1}\sum^2_{\bar s=1}\sum^3_{m=0}
G^{k\bar s}_m\,L_{\boldsymbol\Upsilon_{\!r}}\!(G^{\,m}_{k
\kern 0.5pt\bar s})=2\sum^2_{u=1}\sum^2_{s=1}
L_{\boldsymbol\Upsilon_{\!r}}\!(d_{s\kern 0.3pt u})
\,d^{\kern 0.5pt us}\,+\\
+\,2\sum^2_{\bar u=1}\sum^2_{\bar u=1}
L_{\boldsymbol\Upsilon_{\!r}}\!(\bd_{\bar s\kern 0.3pt \bar u})
\,\bd^{\kern 0.5pt\bar u\bar s}
-\sum^3_{m=0}\sum^3_{n=0}
g^{nm}\,L_{\boldsymbol\Upsilon_{\!r}}\!(\,g_{nm}).
\endgathered
\quad
\mytag{5.43}
$$
Now we substitute \mythetag{5.36} and \mythetag{5.43} into 
\mythetag{5.35}. This yields
$$
\gathered
\sum^2_{k=1}\Alpha^i_{\kern 0.5pt rk}\,G^{k\kern 0.5pt \bar i}_p
+\sum^2_{\bar k=1}\overline{\Alpha^{\bar i}_{r\bar k}}
\,G^{i\kern 0.5pt \bar k}_p
=\sum^3_{q=0}\frac{G^{i\kern 0.5pt\bar i}_q\,\Gamma^q_{rp}}
{2}\,+\\
+\sum^3_{q=0}\frac{G^{i\kern 0.5pt \bar i}_p\,\,\Gamma^q_{rq}}{2}
-\frac{1}{2}L_{\boldsymbol\Upsilon_{\!r}}\!(G^{i\kern 0.5pt\bar i}_p)
-\sum^3_{m=0}\sum^3_{n=0}
\frac{g^{nm}\,L_{\boldsymbol\Upsilon_{\!r}}\!(\,g_{nm})}{4}
\,G^{i\kern 0.5pt \bar i}_p\,-\\
-\sum^2_{k=1}\sum^2_{\bar s=1}\sum^3_{m=0}
\frac{d^{\kern 0.5pt ik}\,\bd^{\kern 0.5pt \bar i\bar s}}{2}
\,\Gamma_{rmp}\,G^{\,m}_{k\kern 0.5pt\bar s}
-\sum^2_{k=1}\sum^2_{\bar s=1}\sum^3_{m=0}
\frac{d^{\kern 0.5pt ik}\,\bd^{\kern 0.5pt \bar i\bar s}
\,g_{p\kern 0.7pt m}}{2}\,L_{\boldsymbol\Upsilon_{\!r}}\!(
G^{\,m}_{k\kern 0.5pt\bar s})\,+\\
+\sum^2_{\bar u=1}\sum^2_{\bar s=1}
\frac{L_{\boldsymbol\Upsilon_{\!r}}\!(\bd_{\bar s\kern 0.3pt\bar u})
\,\bd^{\kern 0.5pt\bar u\bar s}
\,G^{i\kern 0.5pt \bar i}_p}{4}
+\sum^2_{u=1}\sum^2_{s=1}
\frac{L_{\boldsymbol\Upsilon_{\!r}}\!(d_{s\kern 0.3pt u})
\,d^{\kern 0.5pt us}\,G^{i\kern 0.5pt \bar i}_p}{4}.
\endgathered
\qquad
\mytag{5.44}
$$\par
     Now let's study the fifth term in the right hand side of the
equality \mythetag{5.44}. Applying the identity \mythetag{5.40} to
\mythetag{5.44}, we derive
$$
\hskip -2em
\sum^2_{k=1}\sum^2_{\bar s=1}\sum^3_{m=0}
\frac{d^{\kern 0.5pt ik}
\,\bd^{\kern 0.5pt \bar i\bar s}}{2}
\,\Gamma_{rmp}\,G^{\,m}_{k\kern 0.5pt\bar s}
=\sum^3_{q=0}\sum^3_{m=0}\frac{\Gamma_{rmp}}{2}
\,g^{mq}\,G^{i\kern 0.5pt\bar i}_q.
\mytag{5.45}
$$
In order to transform \mythetag{5.45} we use the formula 
\mythetag{5.25}. This yields
$$
\hskip -2em
\gathered
\sum^2_{k=1}\sum^2_{\bar s=1}\sum^3_{m=0}
\frac{d^{\kern 0.5pt ik}
\,\bd^{\kern 0.5pt \bar i\bar s}}{2}
\,\Gamma_{rmp}\,G^{\,m}_{k\kern 0.5pt\bar s}
=\\
=\sum^3_{q=0}\sum^3_{m=0}\frac{g^{qm}\,
L_{\boldsymbol\Upsilon_{\!r}}\!(g_{p\kern 0.5pt m})}{2}
\,G^{i\kern 0.5pt\bar i}_q
-\sum^3_{q=0}\frac{\Gamma^q_{rp}}{2}
\,G^{i\kern 0.5pt\bar i}_q.
\endgathered
\mytag{5.46}
$$
Substituting \mythetag{5.46} back into the formula \mythetag{5.44},
we get the following equality:
$$
\gathered
\sum^2_{k=1}\Alpha^i_{\kern 0.5pt rk}\,G^{k\kern 0.5pt \bar i}_p
+\sum^2_{\bar k=1}\overline{\Alpha^{\bar i}_{r\bar k}}
\,G^{i\kern 0.5pt \bar k}_p
=\sum^3_{q=0}G^{i\kern 0.5pt\bar i}_q\,\Gamma^q_{rp}\,+\\
+\sum^3_{q=0}\frac{G^{i\kern 0.5pt \bar i}_p\,\,\Gamma^q_{rq}}{2}
-\frac{1}{2}L_{\boldsymbol\Upsilon_{\!r}}\!(G^{i\kern 0.5pt\bar i}_p)
-\sum^3_{m=0}\sum^3_{n=0}
\frac{g^{nm}\,L_{\boldsymbol\Upsilon_{\!r}}\!(\,g_{nm})}{4}
\,G^{i\kern 0.5pt \bar i}_p\,-\\
-\sum^3_{q=0}\sum^3_{m=0}\frac{g^{qm}\,
L_{\boldsymbol\Upsilon_{\!r}}\!(g_{p\kern 0.5pt m})}{2}
\,G^{i\kern 0.5pt\bar i}_q
-\sum^2_{k=1}\sum^2_{\bar s=1}\sum^3_{m=0}
\frac{d^{\kern 0.5pt ik}\,\bd^{\kern 0.5pt \bar i\bar s}
\,g_{p\kern 0.7pt m}}{2}\,L_{\boldsymbol\Upsilon_{\!r}}\!(
G^{\,m}_{k\kern 0.5pt\bar s})\,+\\
+\sum^2_{\bar u=1}\sum^2_{\bar s=1}
\frac{L_{\boldsymbol\Upsilon_{\!r}}\!(\bd_{\bar s\kern 0.3pt\bar u})
\,\bd^{\kern 0.5pt\bar u\bar s}
\,G^{i\kern 0.5pt \bar i}_p}{4}
+\sum^2_{u=1}\sum^2_{s=1}
\frac{L_{\boldsymbol\Upsilon_{\!r}}\!(d_{s\kern 0.3pt u})
\,d^{\kern 0.5pt us}\,G^{i\kern 0.5pt \bar i}_p}{4}.
\endgathered
\qquad
\mytag{5.47}
$$\par 
     Let's proceed to the second term in the right hand side of 
\mythetag{5.47}. From \mythetag{5.4} due to the skew symmetry
$c^{\,k}_{ij}=-c^{\,k}_{j\kern 0.5pt i}$ we derive the formula
$$
\hskip -2em
\sum^3_{q=0}\Gamma^q_{rq}=\sum^3_{q=0}\sum^3_{m=0}
\frac{g^{\kern 0.5pt q\kern 0.5pt m}
\,L_{\boldsymbol\Upsilon_{\!r}}\!(g_{q\kern 0.5pt m})}{2}.
\mytag{5.48}
$$
Substituting \mythetag{5.48} back into \mythetag{5.47}, we find 
that the second term in the right hand side of this formula
cancels the fourth term over there. So we have 
$$
\gathered
\sum^2_{k=1}\Alpha^i_{\kern 0.5pt rk}\,G^{k\kern 0.5pt \bar i}_p
+\sum^2_{\bar k=1}\overline{\Alpha^{\bar i}_{r\bar k}}
\,G^{i\kern 0.5pt \bar k}_p
=\sum^3_{q=0}G^{i\kern 0.5pt\bar i}_q\,\Gamma^q_{rp}
-\frac{1}{2}L_{\boldsymbol\Upsilon_{\!r}}\!(G^{i\kern 0.5pt
\bar i}_p)\,-\\
-\sum^3_{q=0}\sum^3_{m=0}\frac{g^{qm}\,
L_{\boldsymbol\Upsilon_{\!r}}\!(g_{p\kern 0.5pt m})}{2}
\,G^{i\kern 0.5pt\bar i}_q
-\sum^2_{k=1}\sum^2_{\bar s=1}\sum^3_{m=0}
\frac{d^{\kern 0.5pt ik}\,\bd^{\kern 0.5pt \bar i\bar s}
\,g_{p\kern 0.7pt m}}{2}\,L_{\boldsymbol\Upsilon_{\!r}}\!(
G^{\,m}_{k\kern 0.5pt\bar s})\,+\\
+\sum^2_{\bar u=1}\sum^2_{\bar s=1}
\frac{L_{\boldsymbol\Upsilon_{\!r}}\!(\bd_{\bar s\kern 0.3pt\bar u})
\,\bd^{\kern 0.5pt\bar u\kern 0.3pt\bar s}
\,G^{i\kern 0.5pt \bar i}_p}{4}
+\sum^2_{u=1}\sum^2_{s=1}
\frac{L_{\boldsymbol\Upsilon_{\!r}}\!(d_{s\kern 0.3pt u})
\,d^{\kern 0.5pt u\kern 0.3pts}\,G^{i\kern 0.5pt \bar i}_p}{4}.
\endgathered
\qquad
\mytag{5.49}
$$
The next step is to transform the fourth term in the right hand side
of the formula \mythetag{5.49}. For this purpose we recall that
$L_{\boldsymbol\Upsilon_{\!r}}$ is the differential operator
\mythetag{2.15}. We apply the Leibniz rule for this operator. As a 
result we get
$$
\gather
\sum^2_{k=1}\sum^2_{\bar s=1}\sum^3_{m=0}
d^{\kern 0.5pt ik}\,\bd^{\kern 0.5pt \bar i\bar s}
\,g_{p\kern 0.7pt m}\,L_{\boldsymbol\Upsilon_{\!r}}\!(
G^{\,m}_{k\kern 0.5pt\bar s})
=\sum^2_{k=1}\sum^2_{\bar s=1}\sum^3_{m=0}
L_{\boldsymbol\Upsilon_{\!r}}\!(d^{\kern 0.5pt ik}
\,\bd^{\kern 0.5pt \bar i\bar s}\,g_{p\kern 0.7pt m}
\,G^{\,m}_{k\kern 0.5pt\bar s})\,-\\
-\sum^2_{k=1}\sum^2_{\bar s=1}\sum^3_{m=0}
L_{\boldsymbol\Upsilon_{\!r}}\!(d^{\kern 0.5pt ik})
\,\bd^{\kern 0.5pt \bar i\bar s}\,g_{p\kern 0.7pt m}
\,G^{\,m}_{k\kern 0.5pt\bar s}
-\sum^2_{k=1}\sum^2_{\bar s=1}\sum^3_{m=0}
d^{\kern 0.5pt ik}
\,L_{\boldsymbol\Upsilon_{\!r}}\!(\bd^{\kern 0.5pt
\bar i\bar s})\,g_{p\kern 0.7pt m}
\,G^{\,m}_{k\kern 0.5pt\bar s}\,-\\
-\sum^2_{k=1}\sum^2_{\bar s=1}\sum^3_{m=0}
d^{\kern 0.5pt ik}
\,\bd^{\kern 0.5pt\bar i\bar s}\,
L_{\boldsymbol\Upsilon_{\!r}}\!(g_{p\kern 0.7pt m})
\,G^{\,m}_{k\kern 0.5pt\bar s}
=L_{\boldsymbol\Upsilon_{\!r}}\!(G^{i\kern 0.5pt\bar i}_p)
-\sum^2_{k=1}\sum^2_{q=1}
L_{\boldsymbol\Upsilon_{\!r}}\!(d^{\kern 0.5pt ik})\,
d_{kq}\,\times\\
\times\,G^{q\kern 0.5pt\bar i}_p-\sum^2_{\bar s=1}
\sum^2_{\bar r=1}L_{\boldsymbol\Upsilon_{\!r}}\!(
\bd^{\kern 0.5pt\bar i\kern 0.5pt\bar s})\,
\bd_{\bar s\bar r}\,G^{i\kern 0.5pt\bar r}_p
-\sum^3_{m=0}\sum^3_{q=0}
L_{\boldsymbol\Upsilon_{\!r}}\!(g_{p\kern 0.7pt m})
\,g^{qm}\,G^{i\kern 0.5pt\bar i}_q.
\endgather
$$
Note that the quantities $L_{\boldsymbol\Upsilon_{\!r}}\!(
d^{\kern 0.5pt ik})$ and $L_{\boldsymbol\Upsilon_{\!r}}\!(
\bd^{\kern 0.5pt\bar i\kern 0.5pt\bar s})$ form two skew-symmetric 
$2\times 2$ matrices. It is known that any skew-symmetric $2\times 2$
matrix is proportional to any other nonzero skew-symmetric $2\times 2$ 
matrix. In particular, we can write 
$$
\xalignat 2
\hskip -2em
L_{\boldsymbol\Upsilon_{\!r}}\!(
d^{\kern 0.5pt ik})=U_r\,d^{\kern 0.5pt ik},
&&L_{\boldsymbol\Upsilon_{\!r}}\!(
\bd^{\kern 0.5pt\bar i\kern 0.5pt\bar s})
=\bar U_r\,\bd^{\kern 0.5pt\bar i\kern 0.5pt\bar s}.
\mytag{5.50}
\endxalignat
$$
The coefficients $U_r$ and $\bar U_r$ in \mythetag{5.50} can be easily 
calculated:
$$
\hskip -2em
\aligned
&U_r=\sum^2_{i=1}\sum^2_{k=1}
\frac{L_{\boldsymbol\Upsilon_{\!r}}\!(d^{\kern 0.5pt i\kern 0.5ptk})}
{2}\,d_{k\kern 0.5pt i}=-\sum^2_{u=1}\sum^2_{s=1}
\frac{L_{\boldsymbol\Upsilon_{\!r}}\!(d_{s\kern 0.3pt u})
\,d^{\kern 0.5pt u\kern 0.3pt s}}{2},\\
&\bar U_r=\sum^2_{\bar i=1}\sum^2_{\bar s=1}
\frac{L_{\boldsymbol\Upsilon_{\!r}}\!(\bd^{\kern 0.5pt\bar i\kern 0.5pt\bar s})}
{2}\,\bd_{\bar s\bar i}
=-\sum^2_{\bar u=1}\sum^2_{\bar s=1}
\frac{L_{\boldsymbol\Upsilon_{\!r}}\!(\bd_{\bar s\kern 0.3pt\bar u}
)\,\bd^{\kern 0.5pt\bar u\kern 0.3pt\bar s}}{2}.
\endaligned
\mytag{5.51}
$$
Substituting \mythetag{5.51} back into \mythetag{5.50} and using 
\mythetag{5.50}, we obtain 
$$
\gathered
\sum^2_{k=1}\sum^2_{\bar s=1}\sum^3_{m=0}
d^{\kern 0.5pt ik}\,\bd^{\kern 0.5pt \bar i\bar s}
\,g_{p\kern 0.7pt m}\,L_{\boldsymbol\Upsilon_{\!r}}\!(
G^{\,m}_{k\kern 0.5pt\bar s})
=L_{\boldsymbol\Upsilon_{\!r}}\!(G^{i\kern 0.5pt\bar i}_p)
+\sum^2_{u=1}\sum^2_{s=1}
\frac{L_{\boldsymbol\Upsilon_{\!r}}\!(d_{s\kern 0.3pt u})
\,d^{\kern 0.5pt u\kern 0.3pt s}}{2}\,\times\\
\times\,G^{i\kern 0.5pt\bar i}_p
+\sum^2_{\bar u=1}\sum^2_{\bar s=1}
\frac{L_{\boldsymbol\Upsilon_{\!r}}\!(\bd_{\bar s\kern 0.3pt\bar u}
)\,\bd^{\kern 0.5pt\bar u\kern 0.3pt\bar s}}{2}\,G^{i\kern 0.5pt
\bar i}_p-\sum^3_{m=0}\sum^3_{q=0}
g^{qm}\,L_{\boldsymbol\Upsilon_{\!r}}\!(g_{p\kern 0.7pt m})
\,G^{i\kern 0.5pt\bar i}_q.
\endgathered
\mytag{5.52}
$$
The last step is to substitute \mythetag{5.52} back into \mythetag{5.49}. 
As a result we get
$$
\hskip -2em
\sum^2_{k=1}\Alpha^i_{\kern 0.5pt rk}\,G^{k\kern 0.5pt \bar i}_p
+\sum^2_{\bar k=1}\overline{\Alpha^{\bar i}_{r\bar k}}
\,G^{i\kern 0.5pt \bar k}_p
=\sum^3_{q=0}G^{i\kern 0.5pt\bar i}_q\,\Gamma^q_{rp}
-L_{\boldsymbol\Upsilon_{\!r}}\!(G^{i\kern 0.5pt\bar i}_p).
\mytag{5.53}
$$\par 
     Comparing \mythetag{5.53} with \mythetag{5.9}, we see that these 
two formulas are equivalent to each other. This means that for the spinor 
connection with the components \mythetag{1.3} the second concordance 
condition \mythetag{5.7} is also fulfilled. \pagebreak As compared to the 
paper \mycite{1}, in this paper we have verified the formula \mythetag{1.3} 
in some other way, though the calculations here are not less huge than 
those in \mycite{1}.

\enddocument
\end